\documentclass[a4paper,11pt]{amsart}

\usepackage{cellspace}
\usepackage[marginparsep=0.6cm,marginparwidth=3cm]{geometry}

\usepackage{amsmath,amsthm,amsfonts,amssymb,mathtools,mathrsfs,bm}
\usepackage[utf8]{inputenc}

\usepackage[T1]{fontenc}
\usepackage{lmodern}

\usepackage{accents}
\usepackage{nicefrac}
\usepackage{diagbox}
\usepackage{makecell}

\usepackage{tikz-cd}
\usetikzlibrary{babel}
\usetikzlibrary{decorations.pathmorphing}
\usepackage[shortlabels]{enumitem}
\usepackage{float}
\usepackage{subcaption}
\usepackage{graphicx}

\usepackage{xcolor}
\definecolor{colordelink}{rgb}{0,0,0.50}
\definecolor{colordecite}{rgb}{0,0.5,0}
\definecolor{colordeurl}{rgb}{0,0.41,0.5}
\usepackage[pagebackref=true]{hyperref}
\hypersetup{
    colorlinks=true,
    linkcolor=colordelink,
    citecolor=colordecite,
		urlcolor=colordeurl,
}
\usepackage[capitalize,nameinlink,noabbrev]{cleveref}
\usepackage[all]{hypcap}
\usepackage{todonotes}

\usepackage{adjustbox}

\def\im{\operatorname{im}}
\def\dim{\operatorname{dim}}
\def\deg{\operatorname{deg}}

\def\coker{\operatorname{coker}}
\def\codim{\operatorname{codim}}

\def\id{\operatorname{id}}

\def\sign{\operatorname{sgn}}
\def\sgn{\operatorname{sgn}}

\def\tr{\operatorname{tr}}

\newcommand{\Alt}{\textnormal{Alt}}

\newcommand{\CC}{\mathbb{C}}

\newcommand{\ZZ}{\mathbb{Z}}

\newcommand{\eqA}{\mathscr{A}}

\newcommand{\Scal}{\mathcal{S}}

\newcommand{\zetabm}{\bm{\zeta}}

\newcommand{\CCzero}[1]{(\mathbb{C}^{#1},0)}

\newcommand{\CCS}[1]{(\mathbb{C}^{#1},S)}
\newcommand{\GS}[2]{(\mathbb{C}^{#1},S)\rightarrow(\mathbb{C}^{#2},0)}
\newcommand{\Gzero}[2]{(\mathbb{C}^{#1},0)\rightarrow(\mathbb{C}^{#2},0)}

\newcommand{\mdash}{\nobreakdash-\hspace{0pt}}

\newcommand{\barbelow}[1]{\underaccent{\bar}{#1}}

\makeatletter
\newcommand{\tpitchfork}{%
  \vbox{
    \baselineskip\z@skip
    \lineskip-.52ex
    \lineskiplimit\maxdimen
    \m@th
    \ialign{##\crcr\hidewidth\smash{$-$}\hidewidth\crcr$\pitchfork$\crcr}
  }%
}
\makeatother

\theoremstyle{plain}
\newtheorem{theorem}{Theorem}[section]
\newtheorem{lemma}[theorem]{Lemma}
\newtheorem{corollary}[theorem]{Corollary}
\newtheorem{proposition}[theorem]{Proposition}

\newtheorem*{theorem**}{Theorem\theoremnum}
\newenvironment{theorem*}[1][]{%
  \edef\theoremnum{\if\relax\detokenize{#1}\relax\else~#1\fi}
  \begin{theorem**}
}{%
  \end{theorem**}
} 
\newtheorem*{conjecture**}{Conjecture\theoremnum}
\newenvironment{conjecture*}[1][]{%
  \edef\theoremnum{\if\relax\detokenize{#1}\relax\else~#1\fi}
  \begin{conjecture**}
}{%
  \end{conjecture**}
}

\theoremstyle{definition}
\newtheorem{definition}[theorem]{Definition}
\newtheorem{proposition/definition}[theorem]{Proposition/Definition}
\newtheorem{conjecture}[theorem]{Conjecture}
\newtheorem{example}[theorem]{Example}
\newtheorem{oproblem}[theorem]{Open Problem}
\newtheorem{remark}[theorem]{Remark}
\newtheorem{notation}[theorem]{Notation}
\Crefname{notation}{Notation}{Notations}

\theoremstyle{remark}
\newtheorem*{note}{Note}

\begin{document}




\author{R. Giménez Conejero}

\title[Isotypical components of the homology of ICIS and images of map germs]
{Isotypical components of the homology of ICIS and images of deformations of map germs}

\address{
Mid Sweden University, Sidsjövägen 2, 852 33 Sundsvall,
Sweden
}
\email{Roberto.Gimenez@uv.es}


\subjclass[2000]{Primary 58K15; Secondary 58K40, 58K65} \keywords{Deformations of map germs, homology, excellent unfoldings, coalescence, isotypical components}

\begin{abstract} 
We give a simple way to study the isotypical components of the homology of simplicial complexes with actions of finite groups, and use it for Milnor fibers of \textsc{icis}. We study the homology of images of mappings $f_t$ that arise as deformations of complex map germs $f:(\CC^n,S)\to(\CC^p,0)$, with $n<p$, and the behaviour of singularities (instabilities) in this context. We study two generalizations of the notion of image Milnor number $\mu_I$ given by Mond and give a workable way of compute them, in corank one, with Milnor numbers of \textsc{icis}. We also study two unexpected traits when $p>n+1$: stable perturbations with contractible image and homology of $\text{im} f_t$ in unexpected dimensions. We show that Houston's conjecture, $\mu_I$ constant in a family implies excellency in Gaffney's sense, is false, but we give a correct modification of the statement of the conjecture which we also prove.
\end{abstract}

\maketitle
\section{Introduction} 

The context of this work is known as Thom-Mather theory: we deal with homology of images of mappings $f_t$ that arise as deformations of complex map germs $f:(\CC^n,S)\to(\CC^p,0)$, with $n<p$, and the behaviour of singularities (instabilities) under deformations. The most studied case is $p=n+1$, where there are several long-standing open problems. In contrast, very little is known for $p>n$, in general (when the deformation $f_s$ is locally stable, the homology of $\im f_s$ is known, see \cite{Houston1997,Goryunov1993,Goryunov1995}).

We refer to the book \cite{Mond2020} for basic definitions and properties about singularities of mappings: instabilities, stability, $\eqA$-finiteness, versal unfoldings, etc.
\newline

In \cref{sec: action general} we develop the theory of actions of finite groups on a finite simplicial complex $M$ with a \textsl{good} action.  We give several statements regarding the isotypical components of the homology $H_*(M)$ and use it to study the case of  \textsc{icis} (isolated complete intersection singularities) and their Milnor fiber.

This is applied in \cref{sec:symmetry,sec:disgusting,sec:muigeneral} to express the homology of $\im f_t$ in simple terms, $f_t$ being \textit{any}
 deformation of a germ $f$. We do this by means of a well known spectral sequence, the Image-Computing Spectral Sequence (\textsc{icss}), whose input is given by the \textsl{multiple point spaces}, $D^k(\bullet)$, introduced in \cref{sec:disgusting}. $D^k(f)$ are \textsc{icis} if the germ $f$ has corank one, and we use Milnor numbers of a larger set of \textsc{icis} to express $H_*(\im f_t)$ in this case.

The deformations that are locally stable stand out among every other deformation, they play the role of Milnor fibers (generic fibers). If we refer only to (locally) stable deformations, it is customary to use $f_s$ instead of $f_t$. In \cref{sec:disgusting} we study unstable germs such that the stable deformation $f_s$ has contractible image
, and characterize the dimensions where such mono-germs exist in corank one. This is similar to the concept of rigid singularities, in our context, and their existence was unknown (apart from trivial cases when $n\ll p$). This also answers negatively a question posed by D. Mond.

The multiple point spaces control the homology of $\im f_t$ by means of the \textsc{icss}, but very little is known about them. In \cref{sec:symmetry} we give deep results about their behaviour (see \cref{cor: sing iff alt,cor: ft tiene homologia}). This is used in \cref{sec:muigeneral} to give a workable way to compute the homology of $\im f_s$, in particular two generalizations of the image Milnor number $\mu_I$ defined by Mond, for $n<p$. This is highly relevant: the homology of $\im f_s$ is difficult to handle in general and we can translate it into a combination of Milnor numbers of \textsc{icis} (in corank one).

In \cref{sec:muigeneral} we also show that there are deformations $f_t$ such that $\im f_t$ has homology in dimensions where the image of the locally stable deformation, $\im f_s$, has trivial homology.
 We show when this \textsl{unexpected homology} appears and we give a \textsl{conservation principle} of the homology of images of perturbations, modulo this unexpected homology.

\textsl{Houston's conjecture on excellent unfoldings} (see \cite{Houston2010}) is very relevant for \textsl{Whitney equisingularity} problems (e.g., \cite{Gaffney1993,GimenezConejero2022b}). In \cref{sec: houston c} we use the previous developments to show that it is false. We modify the conjecture according to these findings and provide a proof in that adjusted
setting. We also show when the converse is true.


Finally, we give some open problems and conjectures to lead the way on the research of this highly unexplored area in \cref{sec: conjectures}.
\newline

\textbf{\emph{Funding:}}
This work was supported by Ministerio de Ciencia, Innovación y Universidades (MICINN) [PID2021-124577NB-I00].
\newline

\textbf{\emph{Acknowledgements:}}
The author thanks Joan Francesc Tent Jorques and Alexander Moretó Quintana for useful indications regarding representation theory. A special mention is to be made about very fruitful conversations between the author and David Mond.

\section{Actions on simplicial complexes}\label{sec: action general}

\subsection{General results}
We use several notions of representation theory, such as characters and irreducible representations. The following result is a key idea in this work, and it involves a system of linear equations where the coefficients are given by the characters of some finite group $G$.

\begin{lemma}\label{lem:charsle}
Let $G$ be a finite group and consider the system of linear equations
$$ \sum_{\tau_i} \chi_{\tau_i}(\sigma) x_{\tau_i}=b_\sigma,$$
one for each element $\sigma\in G$,
with one unknown $x_{\tau_i}$ for each irreducible representation $\tau_i$ of $G$ and $b_\sigma$ the components of a vector $b\in\CC^G$. 
It has a solution if, and only if, $b_{\sigma}=b_{\sigma'}$ for $\sigma$ conjugate to $\sigma'$, and it is  given by
$$ |G|x_{\tau_i}= \sum_{\sigma\in G}\overline{\chi_{\tau_i}(\sigma)}b_\sigma.$$
\end{lemma}
\begin{proof}
Recall that the subspace of vectors $b\in\CC^G$ with $b_\sigma=b_{\sigma'}$ for conjugated $\sigma$ and $\sigma'$ is precisely given by the class functions.
Observe that the irreducible characters $\chi_{\tau_i}$ of $G$ form an orthonormal basis of the space of class functions for the usual inner product
$$ \left\langle \phi| \psi\right\rangle\coloneqq\frac{1}{|G|}\sum_{\sigma\in G}\phi(\sigma)\overline{\psi(\sigma)}.$$

Then, we fix one irreducible representation $\tau'$ and we take the sum of all equations, each one multiplied by the respective $\overline{\chi_{\tau'}(\sigma)}$:
$$\begin{aligned}
 \sum_{\sigma\in G} \overline{\chi_{\tau'}(\sigma)}\sum_{\tau_i} \chi_{\tau_i}(\sigma) x_{\tau_i}&=\sum_{\sigma\in G}\overline{\chi_{\tau'}(\sigma)}b_\sigma\\
\sum_{\tau_i}  \sum_{\sigma\in G} \overline{\chi_{\tau'}(\sigma)}\chi_{\tau_i}(\sigma) x_{\tau_i}&=\sum_{\sigma\in G}\overline{\chi_{\tau'}(\sigma)}b_\sigma\\
\sum_{\tau_i} |G|\left\langle \chi_{\tau_i}|\chi_{\tau'}\right\rangle x_{\tau_i}&=\sum_{\sigma\in G}\overline{\chi_{\tau'}(\sigma)}b_\sigma.
\end{aligned}$$
The lemma follows from there.
\end{proof}

These kind of systems are very relevant in our study because the solution is completely determined by the $b_\sigma$, as the coefficients are given by the character table of $G$. This is one of the main reasons which explain that, as we will see in \cref{thm: conustedesmui}, the \textsl{multiple point spaces} control completely the homology of the image of a \textsl{sufficiently good mapping} $f$.

Another reason is that we use these multiple point spaces in a spectral sequence converging to the homology of their image $\im f$, the so-called \textit{Image-Computing spectral Sequences} or \textsc{icss} (see \cref{thm: general icss}). 

Nonetheless, the difficulty of using these spectral sequences is that we do not use regular homology to compute the limit of the sequence, but an \textsl{alternating isotypical component}.
 In many cases, working with this isotypical component is very hard since one needs to have control over the action of the group on a topological space that can be complicated. Hence, instead of this elusive component, one would like to work with the whole homology (or the Euler-Poincaré characteristic), which is easier to control in favourable circumstances. This is where \cref{lem:charsle} plays a role.

\begin{definition}\label{def:simpliciallygood}
Consider a simplicial complex $M$. We say that an action of a group $G$ on $M$ is \textit{simplicially good} if it maps simplices to simplices and simplices
 fixed as sets are point-wise fixed.
\end{definition}

Let $G$ be a finite group with a simplicially good action on a finite simplicial complex $M$ and let $M^\sigma$ denote the fixed points of the action of $\sigma$ on $M$. Then
\begin{equation}\label{eq:pfeq}
\chi_G(M)(\sigma)=\chi_{Top}(M^\sigma),
\end{equation}
where $\chi_G$ denotes what we call the \textit{(virtual) character of} $G$ \textit{acting on} $M$
 and $\chi_{Top}$ denotes the usual Euler-Poincaré characteristic. More precisely:
\begin{equation*}
		\chi_G(M)(\sigma)\coloneqq \sum_{i} (-1)^i\tr \left(\hspace{-1.5mm}\begin{tikzcd} H_i(M,\CC)\arrow[r,"\sigma_*"]& H_i(M,\CC)\end{tikzcd}\hspace{-1.5mm}\right).
\end{equation*}

\cref{eq:pfeq} links usual homology and the isotypical components of any representation of $G$,
 and it can be proven using standard arguments (see, for example, \cite[Theorem 2.44]{Hatcher2002} and Wall's comments in \cite[p. 172]{Wall1980a}). A detailed account of this can be found in \cite[Section 6.1.1]{Robertothesis}.

It is convenient to use some notation to work with the different irreducible representations of $G$.
\begin{definition}\label{def: tau characteristic}
Consider a simplicially good action of a finite group $G$ on a finite simplicial complex $M$, as in \cref{eq:pfeq}. If $\tau$ is an irreducible representation of $G$, the \textit{$i$-th $\tau$-Betti number} is the number of copies of $\tau$ in $H_i(M,\CC)$, i.e, it is
$$ \beta_i^\tau(M)\coloneqq \frac{\dim_\CC H_i(M,\CC)^\tau}{\deg(\tau)},$$
where $H_i(M,\CC)^\tau$ denotes the $\tau$-isotypical component
 of $H_i(M,\CC)$ and $\deg(\tau)$ the degree of $\tau$. Similarly, the $\tau$\textit{-characteristic of} $M$ is
$$\chi_\tau(M)\coloneqq \sum_{i} (-1)^i \beta_i^\tau(M).$$
\end{definition}

The reason to define $\beta_i^\tau$ normalized by $\deg(\tau)$ instead of defining it as $\dim_\CC H_i(M,\CC)^\tau$ will become apparent in the following proofs.

\begin{proposition}\label{prop: getting chi tau}
In the conditions of \cref{def: tau characteristic},
\begin{equation}
\chi_{\tau}(M)=\frac{1}{|G|}\sum_{\sigma\in G}\overline{\chi_\tau(\sigma)} \chi_{Top}(M^\sigma).
\label{eq: getting chi tau}
\end{equation} 
\end{proposition}
\begin{proof}
For each $\sigma\in G$, we can expand the term $\chi_G(M)(\sigma)$ in \cref{eq:pfeq} to have
$$ \sum_{i} (-1)^i \sum_{\tau_j}\chi_{\tau_j}(\sigma) \beta_i^{\tau_j}(M)= \chi_{Top}(M^\sigma),$$
where $\tau_j$ runs in the irreducible representations of $G$. Rearranging terms,
$$ \sum_{\tau_j} \chi_{\tau_j}(\sigma)\sum_{i} (-1)^i \beta_i^{\tau_j}(M)= \sum_{\tau_j} \chi_{\tau_j}(\sigma)\chi_{\tau_j}(M)=\chi_{Top}(M^\sigma).$$
The statement follows from the previous equation and \cref{lem:charsle}.
\end{proof}

It is useful to state the case when $M$ and all the $M^\sigma$ are non-empty and have non-trivial homology only in dimension zero and, possibly, in another dimension $d$ and $d^\sigma$ (resp.). To simplify the notation, we set $d$ or $d^\sigma$ to zero if $M$ or $M^\sigma$ only have non-trivial homology in dimension zero (resp.).

\begin{proposition}\label{cor:betatauformula}
With the notation of the previous paragraph,
\begin{equation}\label{eq:betatauformula}
(-1)^d\beta_{d}^{\tau}(M)+\zetabm_{d}\beta_0^{\tau}(M)=\frac{1}{|G|}\sum_{\sigma\in G}\overline{\chi_\tau(\sigma)} \big((-1)^{d^\sigma}\beta_{d^\sigma}(M^\sigma)+\zetabm_{d^\sigma}\beta_0(M^\sigma)\big),
\end{equation}
where $\zetabm_a$ is zero when $a=0$ and one when $a\neq 0$.
\end{proposition}
\begin{proof}
This is an immediate reformulation of \cref{prop: getting chi tau} in terms of $\tau$-Betti numbers.
\end{proof}

\begin{figure}
	\centering
		\includegraphics[width=0.75\textwidth]{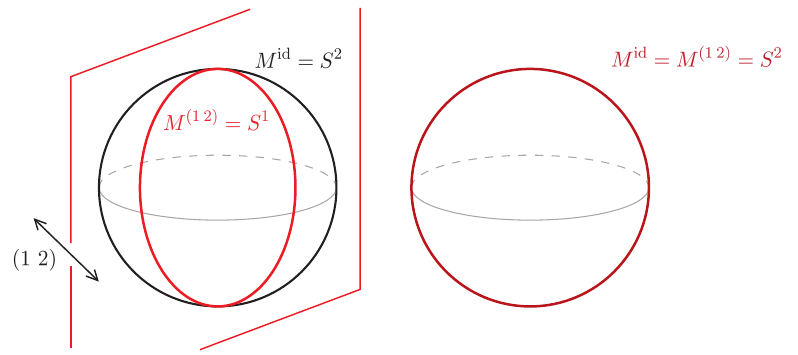}
	\caption{Diagrams of \cref{ex: esfera con accion}.}
	\label{fig:ExampleSphereSigma2}
\end{figure}

\begin{example}\label{ex: esfera con accion}
For the action of $\Sigma_2=\left\{\id, (1\: 2)\right\}$
 on $S^2$ given by a reflection on a plane, the fixed point set of the transposition $(1\: 2)$ are homeomorphic to $S^1$ (see \cref{fig:ExampleSphereSigma2}). Also, the induced action of $\Sigma_2$ on $H_2\left(S^2\right)$ is \textit{alternating}: $(1\: 2)^*\big[S^2\big]=-\big[S^2\big]$. Hence, 
$$\beta_2^\Alt\left(S^2\right)+0= \frac{1}{2}\left(\underbrace{1\cdot(1+1)}_{\id}+\underbrace{(-1)(-1+1)}_{(1\: 2)}\right)=1.$$
In contrast, for the trivial action of $\Sigma_2$,
$$\beta_2^\Alt\left(S^2\right)+0= \frac{1}{2}\left(\underbrace{1\cdot(1+1)}_{\id}+\underbrace{(-1)(1+1)}_{(1\: 2)}\right)=0.$$
\end{example}

\subsection{Complete intersections with isolated singularities}


Now, we study the interesting situation that arises when $M$ has \textsl{vanishing homology}. This happens, for example, when we consider a sufficiently good action on an \textsc{icis} $X$ that is inherited by its Milnor fiber $F_X$. In that case, we are able to compute the \textsl{$\tau$-vanishing homology}, i.e., the $\tau$-Betti number in middle dimension.

A first approach to this situation can be made by assuming that all the $X^\sigma$ are \textsc{icis} with Milnor fiber $F_X^\sigma$, including the cases where $X^\sigma$ are empty. However, as we deal with map germs in the following sections we will have a slightly more general situation that we explain now.  
\newline

In the context of finite map germs $f:(\CC^n,S)\to(\CC^p,0)$ (in particular, $S$ is a finite set), our \textsc{icis} will be a multiple point space $X=D^k(f)$ (see \cref{def: mult spaces f}) with ambient space $\big(\CC^{nk},S^k\big)$. Observe that the \textsc{icis} is a union of germs at several points, so it is rather several \textsc{icis} at several points. In that case, we have an action on
the ambient space by reflections that is inherited by $D^k(f)$, which is given by invariant equations. In particular, the action is also inherited by the Milnor fiber $F_X=D^k(f_s)$.  

The fixed-point spaces $X^\sigma$ of an element $\sigma\in G$ are $D^k(f)^\sigma$. In many cases, they are also \textsc{icis}, since the original \textsc{icis} $X=D^k(f)$ are given by invariant equations and the equations to find the fixed points are those of the hyperplanes where the reflections take place (i.e., the action of the group). In those cases their Milnor fiber are the fixed points of $F_X$ by $\sigma$, i.e., $F_X^\sigma=D^k(f_s)^\sigma$.
However, depending on the dimensions (see \cref{def:expecteddimension}), it can happen that $X^\sigma=D^k(f)^\sigma$ does not have the structure of an \textsc{icis}, but it is still a subset of the set of points $S^k$. In those situations, $F_X^\sigma=D^k(f_s)^\sigma$ is empty. When we allow this to happen, the computations and results need to carry the information about those isolated points, but the main idea of the results is still the same.  
\newline

We extend \cref{def:simpliciallygood} about simplicially good actions to include the setting of multiple point spaces that we will have later and we have explained above. 
\begin{definition}\label{def:sexcellent}
Let $G$ be a finite group and $(X,\Scal)\subset(\CC^N,\Scal)$ an \textsc{icis} on the finite set $\Scal$  (one \textsc{icis} for each point of $\Scal$) with a Milnor fiber $F_X$. We say that $G$ acts in a \textit{simplicially excellent} way if 
\begin{enumerate}[label=(\roman*),font=\itshape]
   \item $G$ has an action on $(\CC^N,\Scal)$ that is inherited by $X$ and by $F_X$; 
	\item the action is simplicially good, i.e., $X$ and $F_X$ admit a simplicial decomposition such that simplices fixed as sets are point-wise fixed; and
   \item\label{it3:sexcellent} for every $\sigma\in G$, $X^\sigma$ are \textsc{icis} with fiber $F_X^\sigma$ or, alternatively, $X^\sigma\subset \Scal$ and $F_X^\sigma=\varnothing$.
   \end{enumerate}
\end{definition}

\begin{note}
The previous definition has the technical difficulty of choosing a Milnor ball such that the action of $G$ on the ambient is inherited to the Milnor fiber $F_X$. In other words, the ball has to be $G$-invariant. This happens in many cases, for example when the action is generated by reflections and $X$ is given by $G$-invariant equations. Indeed, this is the case in the following sections. 
\end{note}

\begin{definition}\label{def:tau mu}
With the conditions of \cref{def:sexcellent}, let $(X,\Scal)$ be an \textsc{icis} of dimension $d$ with a simplicially excellent action of $G$. Consider an irreducible representation $\tau$ of $G$. Then, the \textit{$\tau$-Milnor number} of $X$ is 
$$ \mu^\tau(X;\Scal)\coloneqq (-1)^d \big(\chi_\tau(F_X)-\chi_\tau(X)\big).$$
\end{definition}

\begin{example}
Observe that it is not enough to define $\mu^\tau(X;\Scal)$ as $\beta_d^\tau(F_X)$. The \textsc{icis} $X$ can have dimension zero and, since it is an \textsc{icis} on several points $\Scal$, we need to subtract $\beta_0^\tau(X)$. 
E.g., $K$ copies of the same zero-dimensional \textsc{icis} and $G=\Sigma_K$ the group of permutations of $K$ elements acting by permutation of copies. 
\end{example}

In \cref{it3:sexcellent} of the definition of simplicially excellent actions we include the case when $X^\sigma\subset\Scal$ is only a set of points such that $F^\sigma=\varnothing$. For this reason, we need to extend the definition of Milnor number to that case:
\begin{equation}
\widetilde{\mu}(X^\sigma;\Scal)\coloneqq\begin{cases} 
\mu(X^\sigma;\Scal), &\text{if }X^\sigma\text{ is an \textsc{icis};}\\
-\beta_0(X^\sigma), & \text{otherwise.}
\end{cases}
\label{eq:defmutilde}
\end{equation}

\begin{lemma}\label{lm:computemutau2}
Assume that we have a simplicially excellent action of a finite group $G$ on an \textsc{icis} $(X,\Scal)$.
Then, for any irreducible representation $\tau$,
$$\mu^\tau(X;\Scal)=\frac{1}{|G|}\sum_{\sigma\in G}(-1)^{d-d^\sigma}\overline{\chi_\tau(\sigma)} \widetilde{\mu}(X^\sigma;\Scal),$$
where $d^\sigma$ is the dimension of $X^\sigma$, $d=d^{\id}$ and $\widetilde{\mu}$ is given by \cref{eq:defmutilde}.
\end{lemma}
\begin{proof}
%
Assume that $d>0$. By \cref{cor:betatauformula} for $M=F_X$; and using that $\beta_d^\tau(F_X)=\mu^\tau(X)$,
and that $\beta_0^\tau(F_X)=\beta_0^\tau(X)$ we have
\begin{equation}\label{eq:cambio2eq1}
(-1)^d\mu^\tau(X)+\beta_0^{\tau}(X)=\frac{1}{|G|}\sum_{\sigma\in G}\overline{\chi_\tau(\sigma)} \big((-1)^{d^\sigma}\beta_{d^\sigma}(F_X^\sigma)+\zetabm_{0d^\sigma}\beta_0(F_X^\sigma)\big).
\end{equation}
We can obtain the desired result by eliminating the zero $\tau$-Betti numbers from
 the previous equation, which we can do using \cref{prop: getting chi tau} for $M=X$:
\begin{equation}\label{eq:cambio2eq2}
\beta_0^{\tau}(X)=\frac{1}{|G|}\sum_{\sigma\in G}\overline{\chi_\tau(\sigma)} \beta_0(X^\sigma).
\end{equation}
Subtracting \cref{eq:cambio2eq2} from \cref{eq:cambio2eq1} gives the desired result,
 after taking the notation of $\widetilde{\mu}$ from \cref{eq:defmutilde} and multiplying by $(-1)^d$.

The case $d=0$ follows exactly the same steps, adapting the notation.
\end{proof}

Let $(X,\Scal)$ be an \textsc{icis} with a simplicially excellent action of $G$. Now, consider a $G$-equivariant deformation $X_t$, i.e., $(X_t, G\cdot x)$ is also an \textsc{icis} for any $x\in X_t$ with all the requirements of \cref{def:sexcellent}. Additionally, assume that if $X_0^\sigma$ is not an \textsc{icis} (so it is only points that vanish when one takes the fiber of $X$), neither is the deformation $X_t^\sigma$.
 This is certainly the case for multiple point spaces in the following sections. Then, we have the \textit{conservation principle for the $\tau$-Milnor numbers}. 
\begin{theorem}\label{thm:conservation2}
The $\tau$-Milnor numbers are conserved. 
 In other words, if we have a one-parameter family of \textsc{icis} $X_t$ as in \cref{def:sexcellent}, for $|t|$ small enough,
$$ \mu^\tau(X_0;\Scal)=\begin{cases} \beta_d^\tau(X_t)+\sum_{G x\subseteq X_t} \mu^\tau (X_t;\:G x) & \text{if $d>0$},\\ 
\beta_0^\tau(X_t)+\sum_{G x\subseteq X_t} \mu^\tau (X_t;\:G x)-\beta_0^\tau(X_0)& \text{if $d=0$}.\end{cases}$$ 
\end{theorem}
\begin{proof}
If $d=0$ the result follows immediately from the definitions, so we assume that $d>0$. 

Since we have a simplicially excellent action, we can use the conservation principle of the Milnor number for the \textsc{icis} $X^\sigma$ (see \cite[Theorem 4.2]{Nuno-Ballesteros2018}).  \cref{lm:computemutau2} for $X_0$ and $X_t$ yields 
\begin{align}
	\mu^{\tau}(X_0;\Scal)&=\frac{1}{|G|}\sum_{\sigma\in G}(-1)^{d-d^\sigma}\overline{\chi_\tau(\sigma)} \widetilde{\mu}(X_0^\sigma;\Scal)\text{, and}\label{eq:mux0}\\
	\sum_{G x\subseteq X_t}\mu^{\tau}(X_t;Gx)&=\frac{1}{|G|}\sum_{\sigma\in G}(-1)^{d-d^\sigma}\overline{\chi_\tau(\sigma)} \sum_{G x\subseteq X_t}\widetilde{\mu}(X_t^\sigma; Gx).\label{eq:muxt}
\end{align}
It suffices to prove that subtracting \cref{eq:muxt} from \cref{eq:mux0} gives $\beta_d^\tau(X_t)$. Indeed, by the conservation principle for \textsc{icis}, it gives
$$  \frac{1}{|G|}\left(\sum_{\sigma: d^\sigma>0} (-1)^{d-d^\sigma}\overline{\chi_\tau(\sigma)}\beta_{d^\sigma}(X_t^\sigma) + \sum_{\sigma: d^\sigma = 0} (-1)^{d}\overline{\chi_\tau(\sigma)}\big(\beta_{0}(X_t^\sigma)-\beta_0(X_0^\sigma)\big) \right).$$
This is equal to $\beta_d^\tau(X_t)$. To see it, observe that, in positive dimensions,
\begin{align*}
\chi_\tau(X_0)&=\beta_0^\tau(X_0)=\beta_0^\tau(X_t),\text{ and}\\
\chi_{Top}(X_0^\sigma)&=\beta_0(X_0^\sigma)=\beta_0(X_t^\sigma).
\end{align*}
Hence, using \cref{prop: getting chi tau}, subtracting 
$$ \chi_{\tau}(X_0)=\frac{1}{|G|}\sum_{\sigma\in G}\overline{\chi_\tau(\sigma)} \chi_{Top}(X_0^\sigma) $$
from
$$\chi_{\tau}(X_t)=\frac{1}{|G|}\sum_{\sigma\in G}\overline{\chi_\tau(\sigma)} \chi_{Top}(X_t^\sigma)$$
one obtains the desired equality
$$\beta_d^\tau(X_t)= \frac{1}{|G|}\left(\sum_{\sigma: d^\sigma\neq 0}\overline{\chi_\tau(\sigma)} (-1)^{d-d^\sigma}\beta_{d^\sigma}(X_t^\sigma) + \sum_{\sigma:d^\sigma = 0}\overline{\chi_\tau(\sigma)} (-1)^{d}\big(\beta_{0}(X_t^\sigma)-\beta_0(X_0^\sigma)\big)\right) .$$
 
\end{proof}

\begin{corollary}\label{cor:uppersemicont}
In the situation of \cref{thm:conservation2}, the $\tau$-Milnor numbers are upper semi-continuous: 
$$ \mu^\tau(X_0)\geq \mu^\tau (X_t;x),$$
for all $x\in X_t$ and $|t|$ small enough.
\end{corollary}

The following example shows that the previous theorem does not hold for non-\textsc{icis} and $\widetilde{\mu}$. Also, with the same ideas, one also sees that one needs to assume that if $X_0^\sigma$ is not an \textsc{icis}, neither is the deformation $X_t^\sigma$. This was an observation by the anonymous referee.
\begin{example}
The determinantal singularity $(X,0)\subset(\CC^2,0)$ given by the 2-minors of
\[A_0=\begin{bmatrix}
x^2 & y & 0 \\
0 & x & y 
\end{bmatrix}\]
is not an \textsc{icis}, it has dimension zero and its \textit{Milnor fiber} would be empty (taking  a generic preimage of the equations), so $\widetilde{\mu}(X)=-1$. The deformation $X_t$ given by the 2-minors of
\[A_t=\begin{bmatrix}
x^2-t^2 & y & 0 \\
0 & x-t & y 
\end{bmatrix}\]
consist of the points $p=(t,0)$ and $q=(-t,0)$. It is easy to check that $(X_t,q)$ is an \textsc{icis} but $(X_t,p)$ is not. So, $\widetilde{\mu}(X_t;p)=-1$ and $\widetilde{\mu}(X_t;q)=\mu(X_t;q)=1$ (since the \textsc{icis} has multiplicity two). However,
\[\widetilde{\mu}(X_0)=-1 \neq \beta_0(X_t)+\widetilde{\mu}(X_t;p)+\mu(X_t;q)=2-1+1.\]
\end{example}

\section{Multiple point spaces and strongly contractible germs}\label{sec:disgusting}

In \cref{sec:symmetry}, we use the previous results to study images of stable perturbations of germs $f:(\CC^n,S)\to(\CC^p,0)$, where the germs are $\eqA$-finite and $n<p$. The link between the previous developments and these objects is the Image-Computing Spectral Sequence 
 (\textsc{icss}), where the \textit{multiple point spaces} $D^k(f)$ appear. 

\begin{definition} \label{def: mult spaces f}
The \emph{$k$th-multiple point space} of a mapping or a germ $f$, denoted as $D^k(f)$, is defined as follows:
\begin{itemize}
\item
 Let $f\colon X\rightarrow Y$ be a locally stable mapping between complex manifolds. Then, $D^k(f)$ is equal to the closure of the set of points  $\left(x^{(1)},\dots,x^{(k)}\right)$ in $X^k$ such that $f\left(x^{(i)}\right)=f\left(x^{(j)}\right)$ but $x^{(i)}\neq x^{(j)}$, for all $i\neq j$.

\item When $f\colon (\CC^n,S)\rightarrow(\CC^p,0)$ is a stable germ, $D^k(f)$ is defined analogously but in this case it is a set germ in $\big((\CC^n)^k,S^k\big)$.

\item
Let $f\colon(\CC^n,S)\rightarrow(\CC^p,0)$ be finite, (which implies the existence of a stable unfolding $F(x,u)=\big(f_u(x),u\big)$, see \cite[Proposition 7.2]{Mond2020}). Then, $D^k(f)$ is the complex space germ in $\big((\CC^n)^k,S^k\big)$ given by 
$$
D^k(f)=D^k(F)\cap\left\{u=0\right\}.
$$
\end{itemize}
\end{definition}

\begin{note}
For mono-germs of corank one, $D^k(f)$ can be embedded in $\CC^{n-1}\times\CC^{k}$ (see, for example, \cite[end of p. 371]{Mond1987} or \cite[p. 52]{GimenezConejero2022}).
\end{note}

The theory of \textsc{icss} (Image-Computing Spectral Sequences) starts with a theorem given by Goryunov and Mond in \cite[Proposition 2.3]{Goryunov1993} for rational homology and $\eqA$-finite map germs. Goryunov extended the results on \textsc{icss} with \cite[Corollary 1.2.2]{Goryunov1995} for finite maps and integer homology. Years later, Houston proved \cite[Theorem 5.4]{Houston2007} with a bigger class of maps, an action of an additional group, any coefficients and homology of a pair of spaces. 
 Finally, Cisneros Molina and Mond developed a new simple approach in \cite{CisnerosMolina2022} using double complexes; and Mond and Nuño-Ballesteros gave a self-contained introduction to this topic with the new approach in \cite[Section 10]{Mond2020}. 

The following theorem is stated for \textsl{good maps} in the sense of \cite[Definition 5.1]{Houston2007}, which is a class of maps that contains surjective finite analytic maps between compact subanalytic spaces and also finite and proper simplicial maps between locally-finite simplicial complexes (see \cite[Proposition 5.2 and 5.3]{Houston2007}). This class includes the maps studied in this work.

\begin{theorem}[see {\cite[Theorem 5.4]{Houston2007}}]\label{thm: general icss}
Let $f:X\to Y$ be a continuous map and $\tilde{X}$ a subspace of $X$. Assume that $f:X\to Y$ and $\left.f\right|:\tilde{X}\to Y$ are \textsl{good maps} and that $D^k(\left.f\right|)$ is a subcomplex of $D^k(f)$ for all $k\geq 1$. Then, there exists a spectral sequence
$$ E^1_{p,q}(f,\left.f\right|)=AH_q\big(D^{p+1}(f),D^{p+1}(\left.f\right|);W\big)\Longrightarrow H_{p+q}\big(f(X),f|(\tilde{X});W\big), $$ 
where $W$ is a coefficient group and the differential $d_1$ is induced by the projections $\pi:D^k(f)\to D^{k-1}(f)$ for any $k$.
\end{theorem}

\begin{note}
In practice, one can apply this result to compute the homology of the image of a stable perturbation, $\im f_s$, in two similar ways: with the mapping $f_s$ or considering the pair $\big( \im F, \im f_s\big)$, where $F$ is a versal unfolding and the mapping $\left.f\right|$ in the theorem is $F(\bullet,s)=f_s$.
\end{note}

In this theorem, $AH$ denotes the alternating homology, i.e., the homology of the alternating part of the chain complex: chains so that $\sigma$ acts by $\sgn(\sigma)$. 
 However, when $H_q$ does not have torsion, we can compute the homology and, only then, consider the alternating part (its \textsl{isotypical component}), which is isomorphic to $AH_q$ by \cite[Proposition 10.1]{Mond2020} (cf. \cite[Theorem 2.1.2]{Goryunov1995}).
 This alternating part of the homology is usually denoted as $H^\Alt$. A modern and detailed account of these aspects, multiple point spaces and of the \textsc{icss} can be found in the author's thesis, \cite[Chapter 2]{Robertothesis}. 

Obtaining $AH$, or $H^\Alt$, is very difficult in general because one needs to compute the homology and then the action of the group. 
 Hence, it is desirable to work only with the usual homology of the multiple point spaces, since it is easier to compute and the spaces have good properties in some cases (such as being Milnor fibers of \textsc{icis}, see \cref{thm: Marar-Mond crit} below). 
However, the homology of the multiple point spaces is not enough to determine the rank of the alternating part, as \cref{ex: esfera con accion} illustrates.
 This is discussed in \cref{sec: action general}.
\newline

The multiple point spaces have a very good behaviour in corank one, as the \textit{Marar-Mond criterion} shows. The first version of this criterion was given for mono-germs by Marar and Mond in \cite[Theorem 2.14]{Marar1989}, for multi-germs it was given by Houston in \cite[Theorem 2.4 and Corollary 2.6]{Houston2010}, and for germs of maps with an \textsc{icis} in the source by the author and Nuño-Ballesteros in \cite[Lemma 3.2]{GimenezConejero2022b} (see also \cite[Lemma 4.2.3]{Robertothesis}).

\begin{theorem}\label{thm: Marar-Mond crit}
For a finite germ $f:(\CC^n,S)\rightarrow (\CC^p,0)$ of corank $1$ 
 with $n < p$:
\begin{enumerate}[\itshape(i)]
	\item $f$ is stable if, and only if, $D^k(f)$ is smooth of dimension $p-k(p-n)$, or empty, for any $k\geq1$.\label{item:mararmondi}
	\item \label{it2: Marar-Mond crit}$\eqA_e\textnormal{-codim}(f)$ is finite if, and only if, 
   for each $k$ with $p-k(p-n)\geq 0$, $D^k(f)$ is empty or an \textsc{icis} of dimension $p-k(p-n)$ and if, furthermore, for those $k$ such that $p-k(p-n)< 0$, $D^k(f)$ is a subset of $S^k$, possibly empty.\label{item:mararmondii}
\end{enumerate}
\end{theorem}

As \cref{thm: general icss} shows, the multiple point spaces contain all the information about the image of a mapping by means of the alternating homology. In the case of map germs, this information is summed up in the following objects.

\begin{definition}[cf. {\cite[Definition 3.9]{Houston2010}}]\label{def: mu tau f}
Let $f:\GS{n}{p}$, $n<p$, be $\eqA$-finite of corank 1 and let $F(x,s)=\big(f_s(x),s\big)$ be a stabilisation of $f$. We set the following notation:
\begin{itemize}
\item $s(f)=|S|$, the number of branches of the multi-germ; 
\item $d(f)=\sup\big\{k:\  D^k(f_s)\neq\varnothing\big\}$, where $f_s$ is a stable perturbation of $f$.
\end{itemize}
The \textit{$k$-th $\tau$-Milnor number of $f$}, denoted by $\mu_k^\tau(f)$, is defined as
\[
\mu_k^\tau(f)\coloneqq 
\begin{cases}
\mu^\tau\big(D^k(f)\big),&\  \text{if $k\leq d(f)$} \\[3mm]
	\genfrac(){0pt}{0}{s(f)-1}{d(f)},&\   \text{if $k=d(f)+1$ and $s(f)>d(f)$,}\\[3mm]
	0,&\  \text{otherwise,}
	\end{cases}
\]
following \cref{def:tau mu}.
\end{definition}

Recall that the image of a stable perturbation $f_s$ of an $\eqA$-finite germ $f:(\CC^n,S)\to(\CC^{p},0)$ has the homotopy type of a wedge of spheres. This was proved by Mond for $p=n+1$ in \cite[Theorem 1.4]{Mond1991} and by Houston for any $p>n$ in \cite[Theorem 4.24]{Houston1997}. 
 For $p=n+1$, all the spheres are of dimension $n$ and whose number is the image Milnor number $\mu_I(f)$ \cite[Section 8.3]{Mond2020}. This, and the following proposition, is what motivated the previous definition: 
\begin{proposition}[see {\cite[Definition 3.11]{Houston2010}}]\label{mu}
Let $f\colon(\CC^n,S)\to(\CC^{n+1},0)$ be $\eqA$\mdash finite of corank 1. Then,
	$$\mu_I(f)=\sum_k\mu_k^\Alt(f).$$
\end{proposition}

There is a simple generalization of this result for $p>n$, in general, that is already known (see \cite[Remark 3.12]{Houston2010}). This is studied further in \cref{sec:muigeneral}.

\begin{remark}
The origin of \cref{def: mu tau f} is found in \cite[Definition 3.9]{Houston2010}, in more complicated terms and only for the case $\tau=\Alt$. The first simplification was already made for $\mu_{d(f)+1}^\Alt$ in \cite[Definition 3.2]{GimenezConejero2022} (or \cite[Definition 3.1.14]{Robertothesis}), after observing a property of binomial numbers. For the remaining $\mu_k^\tau$, it is easy to compute the equivalence between both definitions, using
\begin{itemize}
	\item the long exact sequence of the pair $\big(D^k(F), D^k(f_s)\big)$ in $H^\tau_*$;
   \item the fact that $H_0\big(D^k(F),D^k(f_s)\big)$ is also trivial, because there is a point of $D^k(f_s)$ in every connected component of $D^k(F)$; and
   \item the obvious (equivariant) retraction of $D^k(F)$ to $D^k(f)$ for the zero-dimensional case $d_k=0$ (equivalently, $k=\frac{p}{p-n}$).
\end{itemize}
We omit these details since they are a straightforward computation. 
\end{remark}

\begin{remark}\label{rem: s>d}
It is also very important to notice that the inequality $s(f)>d(f)$ can only happen when $d(f)$ reaches the maximum value possible for that pair of dimensions. This was proven in \cite[Lemma 3.3]{GimenezConejero2022}. Observe that this maximum is the number $\kappa\coloneqq\left\lfloor \frac{p}{p-n}\right\rfloor\in\ZZ^+$, by \cref{thm: Marar-Mond crit}. We shall denote this number always with $\kappa$.
\end{remark}

Although the Marar-Mond criterion allows us to give a simple way of computing image Milnor numbers $\mu_I(f)$, it has an unavoidable downside: It could happen that a germ $f$ 
 is unstable but all the spaces $D^k(f)$ are smooth when $p-k(p-n)\geq 0$, forcing to have some $D^k(f)\neq\varnothing$ with $p-k(p-n)<0$. This is an absolutely bad property to have because, for example, any 
 stable perturbation has contractible image in the mono-germ case. This can be seen from the \textsc{icss} of \cref{thm: general icss} for $\im f_s$ (see also \cref{lem: suma mukalt} below). By contrast, for $p=n+1$, it is proven that the image of any stable deformation $\im(f_s)$ of an unstable germ $f$ has non-trivial homology, so it cannot be contractible. This is known as the \textit{weak form of the Mond conjecture} (see \cite[Theorem 3.9]{GimenezConejero2022} and \cite[Theorem 2.15]{GimenezConejero2023a}). 

\begin{figure}
	\centering
		\includegraphics[width=0.80\textwidth]{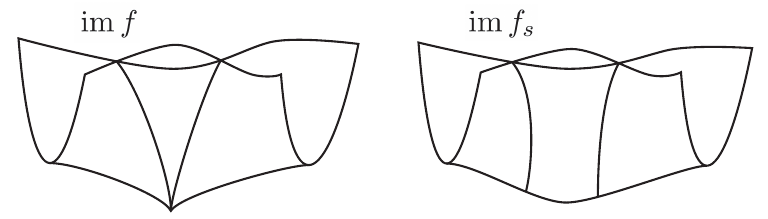}
	\caption{Representation of a strongly contractible instability given by $f$ (left) with its locally stable perturbation $f_s$ (right).}
	\label{fig:strongly contractible}
\end{figure}

\begin{definition}
We will say that a germ $f:(\CC^n,S)\rightarrow (\CC^p,0)$, with $n < p$, is \textit{strongly contractible} if it is unstable, $d(f)\geq s(f)$ and all the spaces $D^k(f)$ are smooth (or, more generally, rigid) of dimension $p-k(p-n)$  
when $p-k(p-n)\geq 0$ (see \cref{fig:strongly contractible,ex:disgusting,rem: ss contr}).
\end{definition}

\begin{remark}
The key point in the previous definition is that the mentioned multiple point spaces $D^k(f)$ must have trivial homology (and, therefore, also alternating homology) after the perturbation corresponding to the stable perturbation $f_s$ of $f$. This is the reason we refer to smoothness or rigidity of such spaces, but it would be enough to say that the spaces $D^k(f_s)$ have trivial alternating homology (cf. \cref{cor: sing iff alt}).
\end{remark}

Observe that germs of corank one which are strongly contractible are always $\eqA$-finite by \cref{it2: Marar-Mond crit} of the Marar-Mond criterion \cref{thm: Marar-Mond crit}.

\begin{notation}\label{not:divdif}
For map germs $f:(\CC^n,0)\to(\CC^p,0)$ of corank one, the equations of the multiple point spaces $D^k(f)$ can be given by the divided differences $\mathit{d}\mathit{d}^j_i$ for $j\leq k$ and $i=1,\dots,p-n+1$ (see \cite[Section 3]{Mond1987} and \cite[Section 9.5]{Mond2020}). If 
$$ f(\barbelow{x},y)=\big(\barbelow{x},f_1(\barbelow{x},y),\dots,f_{p-n+1}(\barbelow{x},y)\big),$$
the divided differences are $\mathit{d}\mathit{d}^1_i(\underline{x},y_1)\coloneqq f_i(\barbelow{x},y_1)$ for $i=1,\dots,p-n+1$, and
$$
\mathit{d}\mathit{d}^k_i(\underline{x},y_1,\dots,y_k)\coloneqq\frac{\mathit{d}\mathit{d}^{k-1}_i(\barbelow{x},y_1,\dots,y_{k-2},y_k)-\mathit{d}\mathit{d}^{k-1}_i(\barbelow{x},y_1,\dots,y_{k-2},y_{k-1})}{y_{k}-y_{k-1}},
$$
 for $i=1,\dots,p-n+1$ and $j\geq 2$. 
\end{notation}

\begin{example}[cf. \cref{ex:n2n-1}]\label{ex:disgusting}
Consider the map germ $f:(\CC^5,0)\to(\CC^8,0)$ given by $$ f(x_1,x_2,x_3,x_4,y)=\big(x_1,x_2,x_3,x_4,\ y^3+x_1y, \ y^4+x_2y,\ y^5+x_3y,\ x_4y+x_1y^2\big).$$

In these dimensions, the unique $k$ such that $p-k(p-n)\geq0$ (the \textsl{expected dimension} of $D^k(f)$) is $k=2$. On the one hand, 
 the equations of $D^2(f)\subset\CC^{5+1}$ are given by the divided differences (see \cref{not:divdif}):
   \begin{align}\label{eq:d2}
\begin{split}
	\mathit{d}\mathit{d}^2_1(\underline{x},y_1,y_2)=& x_1+y_1^2+y_1y_2+y_2^2\\
	\mathit{d}\mathit{d}^2_2(\underline{x},y_1,y_2)=& x_2+y_1^3+y_1^2y_2+y_1y_2^2+y_2^3\\
	\mathit{d}\mathit{d}^2_3(\underline{x},y_1,y_2)=& x_3+y_1^4+y_1^3y_2+y_1^2y_2^2+y_1y_2^3+y_2^4\\
	\mathit{d}\mathit{d}^2_4(\underline{x},y_1,y_2)=& x_4+x_1(y_1+y_2).
	\end{split}
	\end{align}
On the other hand, the equations of $D^3(f)\subset\CC^{5+2}$ are those of $D^2(f)$ and the new divided differences: 
   \begin{align}\label{eq:d3}
\begin{split}
	\mathit{d}\mathit{d}^3_1(\underline{x},y_1,y_2,y_3)=&y_1+y_2+y_3\\
	\mathit{d}\mathit{d}^3_2(\underline{x},y_1,y_2,y_3)=& \sum_{i,j=1,2,3} y_iy_j\\
	\mathit{d}\mathit{d}^3_3\underline{x},y_1,y_2,y_3)=& \sum_{i,j,s=1,2,3} y_iy_jy_s\\
	\mathit{d}\mathit{d}^3_4(\underline{x},y_1,y_2,y_3)=& x_1.
	\end{split}
	\end{align} 
A simple observation of \cref{eq:d2} shows that $D^2(f)$ is smooth, so $f$ is finitely determined by the Marar-Mond criterion, 
 \cref{thm: Marar-Mond crit}. Similarly, it is easy to see from \cref{eq:d2,eq:d3} that $(0,0,0)\in D^3(f)\neq\varnothing$, so $f$ is unstable again by the same criterion. This shows that $f$ is strongly contractible.
\end{example}

\begin{remark}\label{rem: ss contr}
This kind of maps have $\mu_k^\Alt(f)=0$, for every $k$, because the multiple point spaces $D^k(f)$ are smooth, so $D^k(f_s)$ do not have alternating homology. 
 Therefore, the image of the stable perturbation is contractible since it has the homotopy type of a wedge of spheres by \cite[Theorem 4.24]{Houston1997}, 
 and by \cref{lem: suma mukalt} below (which is a simple application 
 of the \textsc{icss}). Furthermore, \cref{cor: sing iff alt} below implies 
 that this is, in fact, a characterization of strongly contractible germs in corank one: they are precisely those that have stable perturbation with contractible image (cf. \cref{fig:strongly contractible}). This is the reason that these maps provide counterexamples to many statements that involve the topology of the stable perturbation, and the reason we have to restrict many results to germs that are not strongly contractible. In fact, we will use map germs similar to the one given in \cref{ex:disgusting} to give counterexamples to several statements.
\end{remark}

We can find some dimensions where there are no corank one germs that are strongly contractible. Indeed, as remarked above, when $p=n+1$ there are no strongly contractible germs (in any corank). To give a better insight of this topic, we introduce now a refinement of the multiple point spaces and the Marar-Mond criterion, which will be important to avoid working with the alternating homology later.

Observe that, in \cref{def: mult spaces f}, there is a straightforward action of $\Sigma_k$ in $D^k(\bullet)$. Indeed, one can give $\Sigma_k$-invariant equations of $D^k(\bullet)$ (see \cite[Section 9.5]{Mond2020} or \cite[Section 2]{Marar1989}) instead of the divided differences used in \cref{ex:disgusting}. This leads to the refinements we were mentioning before for fixed points of these spaces by a permutation $\sigma\in\Sigma_k$. The following definition and theorem have a version for other settings: \cite[Corollary 2.15]{Marar1989} for mono-germs, \cite[Corollary 2.8]{Houston2010} for multi-germs, and \cite[Lemma 3.5]{GimenezConejero2022b} or \cite[Lemma 4.2.10]{Robertothesis} for germs on \textsc{icis}.

\begin{definition}
Let $\Sigma_k$ be the group of permutations of $k$ elements together with its natural action on $D^k(f)$. For any element $\sigma\in\Sigma_k$, we define $D^k(f)^\sigma$ as (the isomorphism type of) the subspace of $D^k(f)$ given by the fixed points of $\sigma$. 
\end{definition}

\begin{remark}\label{rem:decomposition of sigma}
Also, recall that all the elements in $\Sigma_k$ can be decomposed into disjoint cycles in a unique way, called the \textit{cycle shape}. Conversely, if we take a partition of $k$, $\gamma(k)=(r_1,\dots, r_m)$, and $\alpha_i=\texttt{\#}\left\{j:r_j=i\right\}$, one can find an element $\sigma\in\Sigma_k$ such that it can be decomposed into $\alpha_i$ pair-wise disjoint cycles of length $i$ (e.g., the permutation $\sigma=(1\ 2\ 6\ 3)(4\ 5)(8\ 9)(7)(10)$ has cycle shape $(4,2,2,1,1)$,  where $\alpha_4=1$, $\alpha_3=0$, $\alpha_2=2$ and $\alpha_1=2$). 
\end{remark}

\begin{lemma}\label{lem: dkgamma}
With the hypotheses of \cref{thm: Marar-Mond crit} and $\sigma\in\Sigma_k$, we have the following.
\begin{enumerate}[\itshape(i)]
	\item If $f$ is stable, $D^k(f)^\sigma$ is smooth of dimension $p-k(p-n)-k+\sum_i\alpha_i$, or empty.\label{iPD}
	\item $\eqA_e-codim(f)$ is finite if, and only if:\label{iiPD}
	\begin{enumerate}[(a)]
		\item for each $k$ with $p-k(p-n)-k+\sum_i\alpha_i\geq0$, $D^k(f)^\sigma$ is empty or an \textsc{icis} of dimension $p-k(p-n)-k+\sum_i\alpha_i$,
		\item for each $k$ with $p-k(p-n)-k+\sum_i\alpha_i<0$, $D^k(f)^\sigma$ is a subset of $S^k$, possibly empty.
	\end{enumerate}
\end{enumerate}
\label{DP}
\end{lemma}

The following definition follows naturally from the previous result.

\begin{definition}\label{def:expecteddimension}
We will say that the \textit{expected dimension} of $D^k(f)$ (and $D^k(f_t)$) is $d_k\coloneqq p-k(p-n)$, and the \textit{expected dimension} of $D^k(f)^\sigma$ (and $D^k(f_t)^\sigma$) is $d_k^\sigma\coloneqq p-k(p-n)-k+\sum_i\alpha_i$. If we want to emphasize the map, we shall write $d_k(f)$ instead of $d_k$.
\end{definition}

Strongly contractible germs are not very common in corank one when one looks at low dimensions. We restrict to the case of mono-germs for the sake of simplicity, but some of the following results could be generalized to multi-germs using similar ideas. One can see an example of a strongly contractible bi-germ from $\CC^3$ to $\CC^5$ in \cref{ex:multiSS}.
\newline

Let us complete \cref{ex:disgusting} with another example that is stable to show different behaviours. 

\begin{example}\label{ex:n2n-1}
Consider the map germ $f:(\CC^3,0)\to(\CC^5,0)$ given by $$ f(x_1,x_2,y)=\big(x_1,x_2,\ y^3 + x_1y, \ y^4 + x_2y,\ x_2y + y^2\big).$$
In this case $\kappa=2$, $d_2=1$ but $d_2-2+1=0$. The equations of $D^2(f)\subset\CC^{2+2}$ are
\begin{align}\label{eq:d2segundaf}
\begin{split}
	\mathit{d}\mathit{d}^2_1(\underline{x},y_1,y_2)&=x_1 + y_1^2 + y_1 y_2 + y_2^2\\
	\mathit{d}\mathit{d}^2_2(\underline{x},y_1,y_2)&=x_2 + y_1^3 + y_1^2 y_2 + y_1 y_2^2 + y_2^3\\
	\mathit{d}\mathit{d}^2_3(\underline{x},y_1,y_2)&=x_2 + y_1 + y_2.	
	\end{split}
	\end{align}
	Furthermore, one of the equations of $D^3(f)$ is
	$$ \mathit{d}\mathit{d}^3_3(\underline{x},y_1,y_2,y_3)=\frac{\mathit{d}\mathit{d}^2_3(\underline{x},y_1,y_3)-\mathit{d}\mathit{d}^2_3(\underline{x},y_1,y_2)}{y_3-y_2}=1. $$
Hence, $D^3(f)=\varnothing$ and $D^2(f)$ is smooth, so $f$ is stable by Marar-Mond criterion \cref{thm: Marar-Mond crit}.
\end{example}

This example shows that the same ideas of \cref{ex:disgusting} do not work for any pair of dimensions. Another wrong idea would be considering a strongly contractible germ and taking an unfolding, expecting that it is going to be strongly contractible as well. For example, taking $f:(\CC^5,0)\to(\CC^8,0)$ from \cref{ex:disgusting}, every unfolding with one parameter, $F:(\CC^{5+1},0)\to(\CC^{8+1},0)$, is not strongly contractible because $D^2(F)$ keeps being smooth but $D^3(F)$ is always a singular \textsc{icis} of dimension $d_3(F)=0$. 

We present now the characterization of dimensions $(n,p)$ with strongly contractible germs. 

\begin{theorem}\label{thm:disgusting?}
Consider unstable germs $f:(\CC^n,0)\to(\CC^p,0)$ of corank one, $n<p$. Then, there are strongly contractible germs in these dimensions if, and only if,
$$ d_\kappa-\kappa=p- \left\lfloor \frac{p}{p-n}\right\rfloor(p-n+1)\geq0.$$
\end{theorem}
\begin{proof}
The proof is based on a careful inspection of the equations of $D^k(f)$. \cref{ex:disgusting,ex:n2n-1} serve as a good references to follow the argument.
\newline

A strongly contractible germ necessarily has $D^k(f)$ smooth for $k\leq\kappa=\left\lfloor \frac{p}{p-n}\right\rfloor$ and $D^{\kappa+1}(f)\neq\varnothing$. For this to happen, when we look at the equations of the $D^k(f)$ with $k\leq\kappa$ given by divided differences (see \cref{not:divdif}), we need to find independent linear terms in each equation.

On the one hand, remember that a corank one germ can be written as 
$$ f(x_1,\ldots,x_{n-1},y)=\big(x_1,\ldots,x_{n-1},h_{1}(\underline{x},y),\ldots,h_{p-(n-1)}(\underline{x},y)\big).$$
Also, when we study the divided differences that give $D^2(f)$, we have $x_1\ldots,x_{n-1}$ as variables and we introduce $y_1$ and $y_2$. The equations of $D^k(f)$ are those of $D^{k-1}(f)$ and new divided differences, and at each step we add a new variable $y_{k}$. These new divided differences take the form
\[\sum_\delta p_\delta(\underline{x})s_k^\delta(y_1,\dots,y_k),\]
where $p$ are polynomials (possibly constant) and $s_k^\delta$ are the $\Sigma_k$-invariant polynomials
\[s_k^\delta(y_1,\dots,y_k)\coloneqq \sum_{i_1+\dots+i_k=\delta}y_1^{i_1}\cdots y_k^{i_k}.\]
Indeed, for $D^2(f)$ it is obvious, the general case follows from an easy inductive argument on the symmetric part $s_k^\delta$:
\[\frac{s_k^\delta(y_1,\ldots,y_{k-1},y_{k+1})-s_k^\delta(y_1,\ldots,y_{k-1},y_{k})}{y_{k+1}-y_k}=s_{k+1}^{\delta-1}(y_1,\dots,y_{k+1}).\]

On the other hand, if $s_k$ contains a linear term, then it is easy to see that $s_{k+1}$ is a unit. Hence, once we find a linear term on the variables $y_1,\ldots,y_k$ in the equations of $D^k(f)$, $D^{k+1}(f)=\varnothing$. For this reason, once we have found $n-1$ linear terms on $x_1,\ldots,x_{n-1}$ and a linear term on some $y_i$, the next multiple point space is going to be empty. Also, the number of equations necessary to define $D^k(f)$ is $\big(p-(n-1)\big)(k-1)$. 
\newline

Finally, if the number of equations necessary to define $D^\kappa(f)$, $$\big(p-(n-1)\big)\left(\left\lfloor \frac{p}{p-n}\right\rfloor-1\right),$$
is less or equal to 
 $n-1$, we can use the variables $x_1,\ldots,x_{n-1}$ to find a germ so that $D^k(f)$ 
 is smooth for $k\leq\kappa$ and avoids linear terms on $y_i$ except, perhaps, at $D^{\kappa+1}(f)$; so it would be a strongly contractible map. Indeed, if one wants that the linear term $x_i$ appears in a $(k-1)$-iterated divided difference $\mathit{d}\mathit{d}_i^k$ (which define $D^k(\bullet)$) it is enough to include the term $y^{k+1}x_i$ in a component of the map germ (see, for instance, \cref{ex:disgusting}). On the other hand, $y^{k+1}$ gives a linear term in the $(k)$-iterated divided difference. This is enough to give strongly contractible maps whenever
$$ \big(p-(n-1)\big)\left(\left\lfloor \frac{p}{p-n}\right\rfloor-1\right)\leq n-1. $$
This inequality is equivalent to the inequality of the statement.

Furthermore, any strongly contractible germ needs this condition to hold, otherwise $D^{\kappa+1}(f)=\varnothing$ or some $D^k(f)$ is singular for $k\leq\kappa$ (see, for instance, \cref{ex:n2n-1}). This proves the result.
\end{proof}

\begin{remark}
D. Mond posed the following question, that we answer now, regarding a generalization of the Mond conjecture for dimension $(n,p)$ with $p>n+1$: 

\textbf{Q:} \textsl{Can we define a version of the image Milnor number $ \widetilde{\mu_I}$ for any dimensions $(n,p)$ as a linear combination of Betti numbers that is independent of any map germ so that}
$$ \widetilde{\mu_I}(f)\coloneqq \sum_{i>0}a_i\beta_i(\im f_s) \geq \eqA_e-\codim(f) ?$$
 \cref{rem: ss contr,thm:disgusting?} answer negatively to this question, i.e., for a strongly contractible germ $f$ and any linear combination:
$$ \widetilde{\mu_I}(f)\coloneqq \sum_{i>0}a_i\beta_i(\im f_s) =0 < \eqA_e-\codim(f). $$
\end{remark}

\begin{figure}[htb]
\begin{minipage}{0.50\textwidth}
\centering
\includegraphics[width=0.90\textwidth]{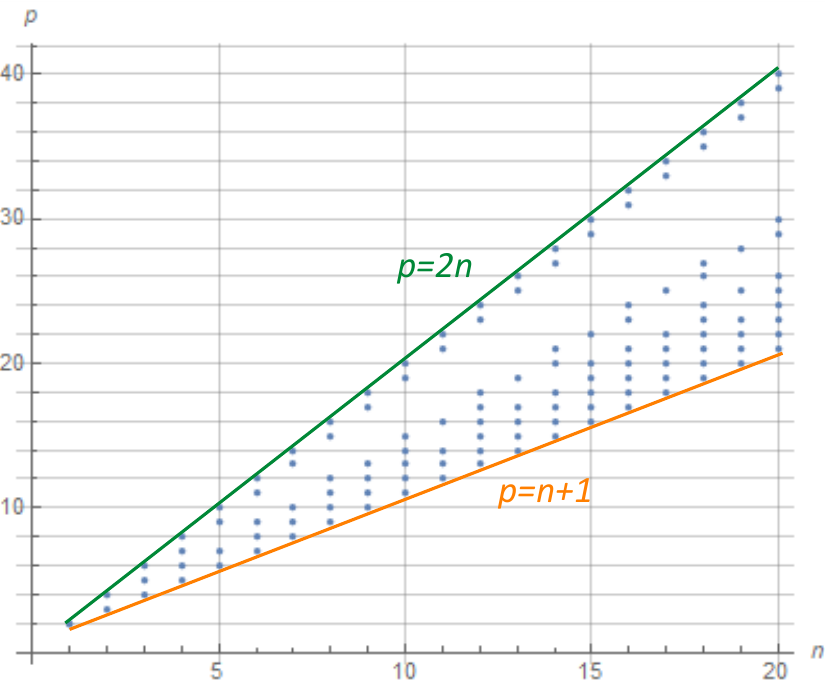}
		\end{minipage}%
		\begin{minipage}{0.50\textwidth}
		\centering
\includegraphics[width=0.90\textwidth]{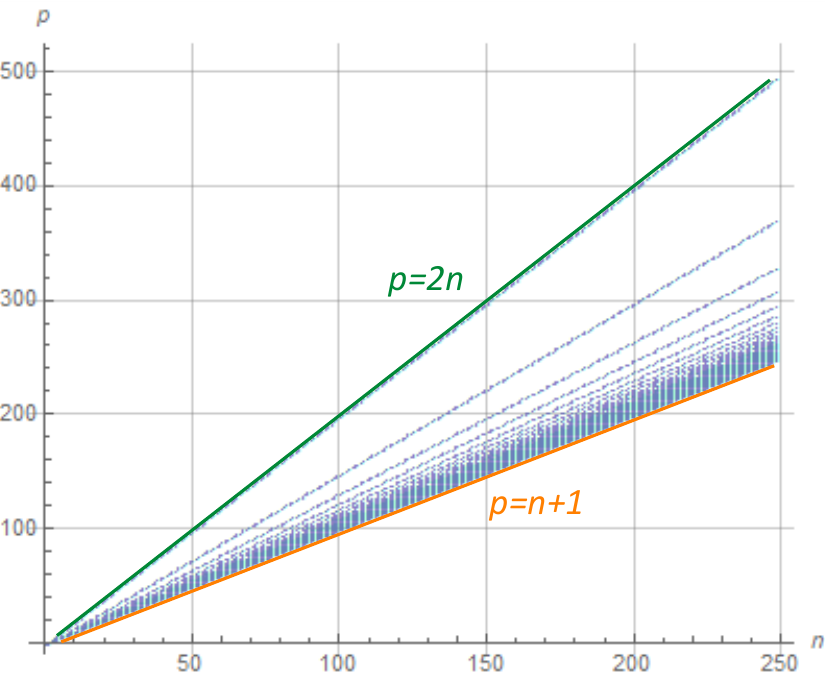}
\end{minipage}
\caption{Plot of the pair of dimensions without strongly contractible germs of corank one, for $n\leq20$ (left) and for $n\leq250$ (right).}
\label{fig:nondisg dimensions}
\end{figure}

\begin{note}\label{rem:disg}
The condition of \cref{thm:disgusting?} implies that a germ $f$ cannot have a perturbation with strongly contractible instabilities if $D^k(f)=\varnothing$ when $d_k<0$. See \cref{fig:nondisg dimensions} for the overall picture of this issue. 
\end{note}

\section{Symmetry of $D^k$}\label{sec:symmetry}

In this section we use the results of \cref{sec: action general} for the multiple point spaces $D^k(f)$ with the natural action of $\Sigma_k$. We are specially interested in the alternating representation $\Alt$, but we can state the conservation of the $\mu^\tau_k(f)$ under deformations of $f$ for any irreducible representation $\tau$ of $\Sigma_k$. 

\begin{theorem}\label{thm: conservation mu tau f}
Let $f:\GS{n}{p}$, $n<p$, be $\eqA$-finite of corank 1. Then, the numbers $\mu^\tau_k(f)$ are conserved for
 $k\leq\frac{p}{p-n}$. More precisely, if $f_t$ is any perturbation of $f$ then, for $|t|$ small enough and $k\leq\frac{p}{p-n}$:
\begin{align*}
\mu^\tau_k(f)=&\beta^\tau_{d_k}\big(D^k(f_t)\big)+\sum_{y\in \im(f_t)}\mu_k^\tau(f_t;y), \text{ if $d_k\neq0$}\\
\mu^\tau_k(f)=&\beta^\tau_{d_k}\big(D^k(f_t)\big)+\sum_{y\in \im(f_t)}\mu_k^\tau(f_t;y)-\beta^\tau_{d_k}\big(D^k(f)\big), \text{ if $d_k=0$};
\end{align*}
 where $\mu_k^\tau(f_t;y)$ is the $k$-th $\tau$-Milnor number of $f$ at $y$.
\end{theorem}
\begin{proof}
The result follows by Marar-Mond criterion (\cref{thm: Marar-Mond crit}), the definition of $\mu^\tau_k$ (\cref{def: mu tau f}) 
 and the conservation principle for the $\tau$-Milnor numbers (\cref{thm:conservation2}). Indeed, we can use the conservation principle since $D^k(f)$ are \textsc{icis} that can be given by $\Sigma_k$-invariant equations (see \cite[Section 9.5]{Mond2020} or \cite[Section 2]{Marar1989}) instead of the divided differences and the action is simplicially excellent, \cref{def:sexcellent}. It is easy to see that the action is simplicially good, but a general result showing this is \cite[Definition 2.8 and Proposition 3.4]{Houston2007} (Houston uses the word \textsl{cellular} instead of \textsl{simplicially good}). 
\end{proof}

The previous result omits the case $ \mu_{d(f)+1}^\tau(f)$, which is non-trivial when $s(f)>d(f)$ (so $k>\frac{p}{p-n}$ by \cref{rem: s>d}). It is false, in general, that this quantity is conservative in a good sense, i.e., that we have something as
$$ \mu_{d(f)+1}^\tau(f)= \mathop{\scalebox{1.5}{\raisebox{-0.2ex}{$\ast$}}} +\sum_{y\in \im(f_t)}\mu_{d(f)+1}^\tau(f_t;y),$$
with the term $\mathop{\scalebox{1.5}{\raisebox{-0.2ex}{$\ast$}}}$ lower semi-continuous. The following example illustrates this fact: $\mu_{d(f)+1}^\Alt$ is neither 
 conservative nor upper semi-continuous.

\begin{example}\label{ex:multiSS}
Consider the following perturbation of a map germ $f:(\CC^3,0)\to(\CC^5,0)$:
$$ f_t(x_1,x_2,y)=(x_1,\ x_2,\ x_1y + y^6 - t^4y^2,\ x_2y + y^5 - t^3y^2,\ x_2y^2 + y^4 - t^2y^2).$$
It is not difficult to see (with a mathematical software, such as \cite{Mathematica}) that $D^2(f_t)\subseteq\CC^{2+2}$
 has singularities at the points
$$ \left(0,\ t^2+2 \sqrt{1-t^2}-2,\ \pm\sqrt{1-\sqrt{1-t^2}},\ \mp\sqrt{1-\sqrt{1-t^2}}\right)$$
and
$$D^3(f_t)=\big\{(0,0;t,0,0),(0,0;0,t,0),(0,0;0,0,t)\big\}.$$
By Marar-Mond criterion, $f$ is $\eqA$-finite. Also, for $t\neq0$, as $f_t(0,0,t)=f_t(0,0,0)=0$, $f_t$ has an instability at the origin given by a bi-germ, which is strongly contractible because $D^2(f_t)$ is smooth at $(0,0;0,t)$, $(0,0;t,0)$ and $(0,0;0,0)$.

We can construct a multi-germ $h=\left\{f,g\right\}:\big(\CC^3,\{0,0'\}\big)\to(\CC^5,0)$ adding an immersion $g$ to the germ $f$ in a way that $h$ is $\eqA$-finite, taking $g$ generic enough. A perturbation $h_t$ can be given by leaving $g$ unchanged and considering the previous perturbation $f_t$ of the non-immersive branch. This perturbation has two unexpected traits:
$$ D^3(h_t)\supset\Sigma_3 (0,0;0,0',t),$$
i.e., it has a free $\Sigma_3$-orbit of points (so $D^3(h_t)$ has alternating homology); and $d(h_0)= s(h_0)=2$ but $2=d(h_t;0)< s(h_t;0)=3$ for $t\neq0$, so $$\mu_{d(h)+1}^{\Alt}(h_0)=0<\mu_{d(h)+1}^{\Alt}(h_t)=1\textnormal{, for } t\neq0.$$ Hence, this is an example of the non-conservation of the $\mu_{d(f)+1}^\tau$ and, in particular, that this number is not upper semi-continuous. This solves a question posed by Houston in \cite[Remarks 4.2 and  6.3]{Houston2010}
\end{example}

Since we are interested in alternated homology and the $k$-th alternating-Milnor numbers $\mu^\Alt_k(f)$ (see also \cref{lem: suma mukalt} below),
 we can use the previous results to compute them. A very trivial case happens when $\Sigma_k$ acts freely on an orbit of a connected component of $D^k(f)$, say $D^k(f)_{free}$ (see \cref{ex: singularidad libre}). In that case, $$\mu^\Alt\big(D^k(f)_{free}\big)=\frac{1}{k!}\mu\big(D^k(f)_{free}\big).$$ If the group does not act freely, things are a bit more complicated and we need to use \cref{lm:computemutau2}.

\begin{figure}
	\centering
		\includegraphics[width=0.25\textwidth]{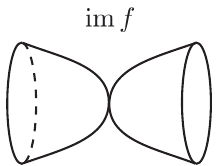}
	\caption{Diagram of the tangent double point $f$ given in \cref{ex: singularidad libre} (cf. \cref{fig:Tangency2}).}
	\label{fig:Tangency1}
\end{figure}

\begin{example}\label{ex: singularidad libre}
It can happen that a germ $f$ has singular $D^k(f)$ but no fixed points for any $\sigma\neq\id$. For example, for $f:\big(\CC^2,\{0,0'\}\big)\to(\CC^3,0)$ given by the \textit{tangent double point}, i.e.,
$$  \left\{\begin{aligned} 
(x,y)&\mapsto (x,y,x^2+y^2)\\
(\widetilde{x},\widetilde{y})&\mapsto\big(\widetilde{x},\widetilde{y}, -(\widetilde{x}^2+\widetilde{y}^2)\big)
\end{aligned}\right. $$
$D^2(f)\subset\CC^{2+2}$ is given by the equations 
$\widetilde{x}=x$,
$\widetilde{y}=y$, and
$x^2+y^2=0$. As the notation $\widetilde{x},\widetilde{y}$ is a simplification to denote multi-germs, there is no fixed point by the permutation $(1\ 2)$ (see \cref{fig:Tangency1}, cf. \cref{fig:Tangency2}).
\end{example}

To tame the sign $(-1)^{d-d^\sigma}$ in \cref{lm:computemutau2}, observe that the difference of dimensions $d_k$ and $d_k^\sigma$ behaves well in the following sense.

\begin{lemma}[see {\cite[Lemma 6.3.2]{Robertothesis}}]\label{lem: parity}
Given an $\eqA$-finite germ $f:\GS{n}{p}$, with $n<p$, $$(-1)^{d_k-d_k^\sigma}=\sgn(\sigma)=\chi_\Alt(\sigma).$$
\end{lemma}
\begin{proof}
First of all, observe that $d_k-d_k^\sigma=k-\sum_i\alpha_i$ provided that $\sigma$ has cycle type $(r_1,\dots, r_m)$ and $\alpha_i=\texttt{\#}\left\{j:r_j=i\right\}$ (see \cref{rem:decomposition of sigma}). 

Assume that the lemma is true for the permutations of $k-1$ elements. To create a permutation of $k$ elements, say $\sigma'$, from a permutation of $k-1$ elements, say $\sigma$, we have to add the new element to any cycle of $\sigma$ or leave that element invariant by $\sigma'$. 

If we add the element to some cycle, the sum $\sum_i \alpha^{\sigma}_i$ transforms into $\sum_i\alpha^{\sigma'}_i$, but notice that we are changing two $\alpha_i$. 
 Therefore, the sum does not change because, for some $i_0$, we have that $\alpha_{i_0}^\sigma =\alpha_{i_0}^{\sigma'}-1 $ and $\alpha_{i_0+1}^\sigma =\alpha_{i_0+1}^{\sigma'}+1 $. Hence, 
$$k-\sum \alpha^{\sigma'}_i=k-\sum \alpha^{\sigma}_i=(k-1)-\sum \alpha^{\sigma}_i+1,$$ 
so the parity of this number changes. Luckily, this operation changes the sign of the permutation: $\sign(\sigma)=-\sign(\sigma')$. 

Finally, if we leave the new element invariant neither the sign nor the difference $k-\sum_i\alpha_i$ changes.

We finish the argument by induction, as it is true for $\Sigma_2$ and $D^2(f)$. 
\end{proof}

\begin{definition}\label{def: mu+-}
The \textit{positive Milnor characteristic} of an \textsc{icis} $(X,\Scal)$ is defined as
$$ \mu^{+0}(X;\Scal)=\mu(X,\Scal)+\beta_0(X),$$
and the \textit{negative Milnor characteristic} as
$$ \mu^{-0}(X;\Scal)=\mu(X;\Scal)-\beta_0(X).$$
If $X$ is 
 empty, both are zero.
\end{definition}

\begin{remark}\label{rem: mu+0-0}
For an \textsc{icis}  $(X,\Scal)$ of dimension zero, $\mu^{+0}(X;\Scal)$ is equal to the number of points of the fiber. Furthermore, $\mu^{-0}(X;\Scal)\geq0$ if each connected component of $X$ has at least one singularity.
\end{remark}

\begin{theorem}\label{conjrgc}
Consider an unstable $\eqA$-finite germ $f:\GS{n}{p}$ of corank one, $n<p$, and a space $D^k(f)\neq\varnothing$ with $d_k\geq0$. Then,
\begin{equation}\label{eqthm: mu +betis}
\mu^\Alt\big(D^k(f)\big)=\frac{1}{k!}\left(\sum_{ d_k^\sigma\geq0}\mu\big(D^k(f)^\sigma\big)-\sum_{ d_k^\sigma<0}(-1)^{d_k^\sigma}\beta_0\big(D^k(f)^\sigma\big)\right),
\end{equation}
where the sums run on the permutations $\sigma\in\Sigma_k$ such that the condition given on $d^\sigma_k$ is satisfied.

If, furthermore, the group $\Sigma_k$ acts freely on the connected components of $D^k(f)$ (i.e., only by permutation of copies of \textsc{icis}), then $\mu^\Alt\big(D^k(f)\big)=\frac{1}{k!}\mu\big(D^k(f)\big)$. Otherwise, we can consider the orbits of those connected components where the group does not act freely and we have the equivalent expression 
\begin{equation}\label{eqthm: mu+-}
\mu^\Alt\big(D^k(f)\big)=\frac{1}{k!}\left(\sum_{d_k^\sigma\geq0\;\textnormal{even}} \mu^{+0}\big( D^k(f)^\sigma\big)+\sum_{d_k^\sigma>0\;\textnormal{odd}} \mu^{-0}\big( D^k(f)^\sigma\big)\right),
\end{equation}
following \cref{def: mu+-}.
\end{theorem}
\begin{proof}
We simply have to rewrite \cref{lm:computemutau2} for $\tau=\Alt$ (we can use it since the action is simplicially excellent, see the proof of \cref{thm: conservation mu tau f}): 
\begin{equation}\label{eq: beware of parity}
\begin{aligned}
\mu^\Alt\big(D^k(f)\big)&= \frac{1}{k!}\sum_{\sigma\in \Sigma_k}(-1)^{d_k-\dim D^k(f)^\sigma}\sign(\sigma) \widetilde{\mu}\big(D^k(f)^\sigma\big)\\
&= \frac{1}{k!}\left(\sum_{d_k^\sigma\geq0}\mu\big(D^k(f)^\sigma\big)-\sum_{d_k^\sigma<0}(-1)^{d_k}\sign(\sigma)\beta_0\big(D^k(f)^\sigma\big)            \right)
\end{aligned}
\end{equation}
where we have simplified signs (using the parity given in \cref{lem: parity}) and we have followed the definition $\widetilde{\mu}\big(D^k(f)^\sigma\big)=-\beta_0\big(D^k(f)^\sigma\big)$ when $d^\sigma_k<0$. 

Observe that, independently of the parity of $d_k$, we have that
\begin{equation}\label{eq: parity is shit}
-\sum_{d_k^\sigma<0}(-1)^{d_k}\sign(\sigma)\beta_0\big(D^k(f)^\sigma\big)= -\sum_{ d_k^\sigma<0 \text{ even }}\beta_0\big(D^k(f)^\sigma\big)+\sum_{ d_k^\sigma<0 \text{ odd }}\beta_0\big(D^k(f)^\sigma\big),
\end{equation}
by \cref{lem: parity}.

\cref{eq: beware of parity,eq: parity is shit} give the first equality of the theorem.

To prove the second one, we only have to show that
\begin{equation}\label{eq:betis0}
 -\sum_{ d_k^\sigma<0}(-1)^{d_k^\sigma}\beta_0\big(D^k(f)^\sigma\big)=\sum_{ d_k^\sigma\geq0}(-1)^{d_k^\sigma}\beta_0\big(D^k(f)^\sigma\big).
\end{equation}
This follows from the fact that there is no alternating homology on $H_*\big(D^k(f)\big)$, provided that the group acts non-freely on connected components. Indeed, $D^k(f)$ retracts to points of the subset $S^k$; but points with some entry repeated at least twice by hypothesis, therefore fixed by some transposition. Hence, there cannot be alternated homology in $H_0$. 
This and \cref{prop: getting chi tau} applied to $D^k(f)$ and $\tau=\Alt$ yield
\begin{equation*}
\begin{aligned}
0 &= \frac{1}{k!}\sum_\sigma \sign(\sigma)\beta_0 \big(D^k(f)^\sigma\big)\\
0&=\sum_\sigma (-1)^{d_k-d_k^\sigma}\beta_0 \big(D^k(f)^\sigma\big)\\
0&=\sum_\sigma (-1)^{d_k^\sigma}\beta_0 \big(D^k(f)^\sigma\big).
\end{aligned}
\end{equation*}
 
\end{proof}

The following theorem sums up most of the behaviour of multiple point spaces and contains a generalization to multi-germs of the results \cite[Theorems 6.3.1 and 6.3.5]{Robertothesis} found in the author's thesis.

\begin{theorem}\label{cor: sing iff alt}
For unstable $\eqA$-finite germ $f:\GS{n}{p}$ of corank one, $n<p$ and spaces $D^k(f)\neq\varnothing$ with $d_k\geq0$, the following are equivalent:
\begin{enumerate}[label=(\roman*),font=\itshape]
	\item $D^k(f)$ has a singularity; \label{item: singularity}
	\item $D^k(f_s)$ has non-trivial homology in middle dimension, $f_s$ stable; \label{item: homology}
	\item $D^{k+j}(f)$ has a singularity if it is not empty, for every $j\geq0$ such that $d_{k+j}\geq0$; \label{item: singularity k+s}
	\item every space $D^k(f)^\sigma$ has a singularity if it is not empty and $d_k^\sigma\geq0$; \label{item: dksigma sing todos}
	\item some space $D^k(f)^\sigma$ has a singularity, where $d_k^\sigma\geq0$;  \label{item: dksigma sing}
	\item $D^k(f_s)$ has non-trivial alternating homology in middle dimension. \label{item: alt homology}
\end{enumerate}
Furthermore, if $d_k^\sigma<0$ and $D^k(f)^\sigma$ is not empty, we have all of the above.
\end{theorem}
\begin{proof}
We can exclude the orbits  of those connected components where the group acts freely, i.e., only by permutation of copies of \textsc{icis}, as the result is trivial there. 
Hence, we can assume that $D^k(f)$ is given by a single orbit of a connected component and the group does not act freely on them.

That \cref{item: singularity,item: homology} are equivalent is known, a Milnor fiber of an \textsc{icis} is contractible if, and only if, the \textsc{icis} is smooth. See, for example, the classical result on the Milnor-Tjurina inequality \cite[Theorem]{Looijenga1985}. 

For the equivalence between \cref{item: singularity,item: singularity k+s}, assume that $f$ is a mono-germ. Then, it is given by the fact that the regular sequence defining the \textsc{icis} $D^{k+j}(f)$ contains the regular sequence of the \textsc{icis} $D^{k}(f)$. 
 For multi-germs, also, this relation is true between components of the multiple point spaces taking into account combinations of branches of the multi-germ. For example, if a component $D\subseteq D^k(f)$ is given by points of, say, only one branch, then there is another component $D'\subseteq D^{k+1}(f)$ given by the equations of $D$ and adding the intersection with another branch of $f$.

That \cref{item: singularity} implies \cref{item: dksigma sing todos} is similar: the regular sequence of the \textsc{icis} $D^k(f)^\sigma$ (recall the stronger version of the Marar-Mond criterion \cref{lem: dkgamma}) include that of $D^k(f)$, since we only need to add equations identifying some variables to have $D^k(f)$ from $D^k(f)^\sigma$. The converse is obvious.

\cref{item: singularity} implies \cref{item: dksigma sing} with $\sigma=\id$. For the converse we can use \cref{eq:pfeq}:
$$\chi_{\Sigma_k}\big(D^k(f)\big)(\sigma)=\chi_{Top}(D^k(f)^\sigma).$$
Assume, for the sake of the contradiction, that $D^k(f)$ is smooth but there is some $\sigma$ such that $D^k(f)^\sigma$ is singular, $d_k^\sigma\geq0$. Then, because $D^k(f)$ is smooth, the action on the connected components of $D^k(f)$ and its Milnor fiber $D^k(f_s)$ are the same, i.e., 
$$\chi_{\Sigma_k}\big(D^k(f)\big)(\sigma)=\chi_{\Sigma_k}\big(D^k(f_s)\big)(\sigma),$$
 which is a contradiction with 
$$\chi_{Top}(D^k(f)^\sigma)\neq\chi_{Top}(D^k(f_s)^\sigma),$$
since $D^k(f_s)^\sigma$ is the Milnor fiber of the \textsc{icis} $D^k(f)$ (recall \cref{lem: dkgamma}). 


The equivalence between \cref{item: singularity,item: alt homology} is, now, easy. The converse implication is trivial. The direct implication follows from \cref{item: dksigma sing todos}: if every $D^k(f)^\sigma$ is empty when $d_k^\sigma<0$, the implication follows easily from \cref{eqthm: mu +betis}. Otherwise, assume 
 that there exists $\sigma_0$ such that $D^k(f)^{\sigma_0}$ is non-empty  and $d_k^{\sigma_0}=0$. In that case, the implication follows from \cref{eqthm: mu+-} and \cref{item: dksigma sing todos}, because (recall \cref{rem: mu+0-0}) $$\mu^{\pm0}\big(D^k(f)^\sigma\big)\geq0,$$ 
and $\mu^{+0}\big(D^k(f)^{\sigma_0}\big)>0$. 

Indeed, the assumption always holds if there is a $D^k(f)^{\sigma_-}\neq\varnothing$ so that $d_k^{\sigma_-}<0$. If we consider the permutation $\sigma_-$ and we fix one entry of its cycle decomposition to have a new $\alpha'$, we have another $D^k(f)^{\sigma'}\neq\varnothing$ and $0\leq d_k^{\sigma'}-d_k^{\sigma_-}\leq1$ (recall \cref{rem:decomposition of sigma,lem: dkgamma}). E.g., $(1\ 2\ 6\ 3)(4\ 5)(8\ 9)(7)(10)$ fixing $3$ and $9$ is the permutation $(1\ 2\ 6)(4\ 5)(3)(8)(9)(7)(10)$. Since $d_k=d_k^{\id}\geq0$, if we iterate this process we reach at some step a $\sigma_0$ so that $D^k(f)^{\sigma_0}\neq\varnothing$ and $d_k^{\sigma_0}=0$.  

Finally, regarding the last statement, it follows from  
a similar argument to that of the equivalence between \cref{item: singularity,item: dksigma sing}.
%
%
 
\end{proof}

To complete the previous theorem, we also have this result.

\begin{theorem}\label{cor: ft tiene homologia}
Let $f:\GS{n}{p}$ be an unstable $\eqA$-finite germ of corank one, $n<p$, and $f_t$ a perturbation. If $\beta_{d_k}\big(D^k(f_t)\big)>0$, with $d_k>0$, then $D^k(f_t)$ has alternating homology. Furthermore,  if $\beta_{d_k}\big(D^k(f_t)\big)>0$,
$$ \mu_k^\Alt(f)> \sum_{y\in \im f_t}\mu_k^\Alt(f_t;y\big).$$
\end{theorem}
\begin{proof}
We can assume that the new homology of $D^k(f_t)$ has appeared in an orbit of a connected component where the group does not act freely, otherwise the group acts by permutation of connected components and the result follows from usual homology arguments. We keep denoting as $D^k(f)$ and $D^k(f_t)$ this orbit of the connected component. 

 For the first statement, then, 
at least one term of \cref{eqthm: mu+-} decreases when we consider $f$ and $f_t$: 
$$\mu^{\pm0}\big(D^k(f)\big)>\sum_{x\in D^k(f_t) }\mu^{\pm0}\big(D^k(f_t);x\big) ,$$
by the conservation principle of the (regular) Milnor number for \textsc{icis} (\cite[Theorem 4.2]{Nuno-Ballesteros2018}) and the hypothesis $\beta_{d_k}\big(D^k(f_t)\big)>0$. This implies that $D^k(f_t)$ must have alternating homology by \cref{thm:conservation2}, because when we compute $\mu^\Alt\big(D^k(f)\big)$ with \cref{eqthm: mu+-} it is strictly bigger than the computation with \cref{eqthm: mu+-} of 
$$\sum_{\Sigma_k x\subseteq D^k(f_t)}\mu^\Alt\big(D^k(f_t);\Sigma_k x\big).$$
This last argument also proves the second statement, recall that $\mu^\Alt\big(D^k(f)\big)=\mu_k^\Alt(f)$ and 
$$\sum_{\Sigma_k x\subseteq D^k(f_t)}\mu^\Alt\big(D^k(f_t);\Sigma_k x\big)=\sum_{y\in \im f_t}\mu_k^\Alt(f_t;y\big).$$

\end{proof}

\section{Homology of deformations and unexpected homology}\label{sec:muigeneral}

\subsection{Stable deformations} 
As it has been done several times (see, for example, \cite{Houston2010}), we can use the \textsc{icss} given in \cref{thm: general icss} and the $\mu_k^\Alt$ to compute the image Milnor number of germs $f:\GS{n}{n+1}$ (see \cref{fig: las mejores tablas del universo}). It is not clear what should be the generalization of this object to map germs $f:\GS{n}{p}$, for arbitrary 
 $n<p$. There is not a correct answer to this question but one should consider some definition that allows generalizing properties of germs when $p=n+1$ to the general case $n<p$ only changing $n+1$ by $p$. 

We give, now, two reasonable generalizations. However, we will see later several reasons that make the first definition seemingly better.

\begin{definition}\label{def:muinui}
Consider an $\eqA$-finite germ $f:\GS{n}{p}$, $n<p$, such that $d_2=2n-p\geq0$. We define its \textit{image Milnor number} as  
$$\mu_I(f)\coloneqq\sum_{i>0} \beta_i\big(\im f_s\big)$$
 for any stable perturbation $f_s$. 
Similarly, the \textit{image vanishing characteristic} of $f$ is
$$\nu_I(f)\coloneqq (-1)^{d_2+1}\sum_{i>0} (-1)^{i}\beta_i\big(\im f_s\big),$$
for any stable perturbation $f_s$. If $d_2<0$ (i.e., $p>2n$) we say that $f$ is \textit{degenerate}. In this case, we take as definition $\mu_I(f)=\nu_I(f)=\beta_0\big(\im(f_s)\big)-1$ (compare \cref{fig: las mejores tablas del universo} and \cref{fig: las segundas mejores tablas del universo}).
\end{definition}

The general definition of $\mu_I$ in \cref{def:muinui} is not new. For example, Houston took this definition in \cite{Houston2010}. The second definition can also be seen as the difference $\chi_{Top}(\im f)-\chi_{Top}(\im f_s)$ with a convenient sign (this approach can also be found in, for example, \cite{NOT}).

\begin{lemma}[cf. {\cite[Remark 3.12]{Houston2010}}]\label{lem: suma mukalt}
Given an $\eqA$-finite germ $f:\GS{n}{p}$ of corank one, $n<p\leq2n$, we have that
\begin{align*}
	  \mu_I(f)&=\sum_k \mu_k^\Alt(f),    \\
   \nu_I(f)&=(-1)^{d_2+1}\left[\sum_{k=2}^{ d(f)} (-1)^{d_k+k-1}\mu_k^\Alt(f)\right]+(-1)^{d(f)+d_2}\mu_{d(f)+1}^\Alt(f).
\end{align*} 
In the degenerate case they coincide and are equal to $\mu^\Alt_{d(f)+1}(f)-1$ if $s(f)>1$, and zero otherwise.
\end{lemma}
\begin{proof}
This is a simple application of the \textsc{icss} shown in \cref{thm: general icss} to compute the homology of the pair $\big(\im F ,\im f_s \big)$ and the definition and comments given in \cref{def: mu tau f} (see also \cref{fig: las mejores tablas del universo,fig: las segundas mejores tablas del universo}).
\end{proof}

\begin{figure}[htb]
\hspace{-0.5cm}
\begin{subfigure}{0.4\textwidth}
\scalebox{0.81}{
\begin{tabular}{c | c@{\hspace{1.3\tabcolsep}}  c@{\hspace{1.3\tabcolsep}}  c@{\hspace{1.3\tabcolsep}}   c@{\hspace{1.3\tabcolsep}}   c@{\hspace{1.3\tabcolsep}}   c@{\hspace{1.3\tabcolsep}}   c@{\hspace{1.3\tabcolsep}}   c@{\hspace{1.3\tabcolsep}}}
		$6   $&\textcolor[rgb]{0.75,0.75,0.75}{0}&\textcolor[rgb]{0.75,0.75,0.75}{0}&\textcolor[rgb]{0.75,0.75,0.75}{0}&\textcolor[rgb]{0.75,0.75,0.75}{0}&\textcolor[rgb]{0.75,0.75,0.75}{0}&\textcolor[rgb]{0.75,0.75,0.75}{0}&\textcolor[rgb]{0.75,0.75,0.75}{0}&\textcolor[rgb]{0.75,0.75,0.75}{0}\\[1.5pt]
		$5   $&\textcolor[rgb]{0.75,0.75,0.75}{0}&$		\mu_2^\Alt\hspace{-2pt}	 $&\textcolor[rgb]{0.75,0.75,0.75}{0}&\textcolor[rgb]{0.75,0.75,0.75}{0}&\textcolor[rgb]{0.75,0.75,0.75}{0}&\textcolor[rgb]{0.75,0.75,0.75}{0}&\textcolor[rgb]{0.75,0.75,0.75}{0}&\textcolor[rgb]{0.75,0.75,0.75}{0}\\[1.5pt]
		$4	 $&\textcolor[rgb]{0.75,0.75,0.75}{0}&\textcolor[rgb]{0.75,0.75,0.75}{0}&$  \mu_3^\Alt\hspace{-2pt}     $&\textcolor[rgb]{0.75,0.75,0.75}{0}&\textcolor[rgb]{0.75,0.75,0.75}{0}&\textcolor[rgb]{0.75,0.75,0.75}{0}&\textcolor[rgb]{0.75,0.75,0.75}{0}&\textcolor[rgb]{0.75,0.75,0.75}{0}\\[1.5pt]
		$3	 $&\textcolor[rgb]{0.75,0.75,0.75}{0}&\textcolor[rgb]{0.75,0.75,0.75}{0}&\textcolor[rgb]{0.75,0.75,0.75}{0}&$    \mu_4^\Alt\hspace{-2pt}    $&\textcolor[rgb]{0.75,0.75,0.75}{0}&\textcolor[rgb]{0.75,0.75,0.75}{0}&\textcolor[rgb]{0.75,0.75,0.75}{0}&\textcolor[rgb]{0.75,0.75,0.75}{0}\\[1.5pt]
		$2	 $&\textcolor[rgb]{0.75,0.75,0.75}{0}&\textcolor[rgb]{0.75,0.75,0.75}{0}&\textcolor[rgb]{0.75,0.75,0.75}{0}&\textcolor[rgb]{0.75,0.75,0.75}{0}&$ 		\mu_5^\Alt\hspace{-2pt}	$&\textcolor[rgb]{0.75,0.75,0.75}{0}&\textcolor[rgb]{0.75,0.75,0.75}{0}&\textcolor[rgb]{0.75,0.75,0.75}{0}\\[1.5pt]
		$1	 $&\textcolor[rgb]{0.75,0.75,0.75}{0}&\textcolor[rgb]{0.75,0.75,0.75}{0}&\textcolor[rgb]{0.75,0.75,0.75}{0}&\textcolor[rgb]{0.75,0.75,0.75}{0}&\textcolor[rgb]{0.75,0.75,0.75}{0}&$ 		\mu_6^\Alt\hspace{-2pt}	$&\textcolor[rgb]{0.75,0.75,0.75}{0}&\textcolor[rgb]{0.75,0.75,0.75}{0}\\[1.5pt] 
		$0	 $&\textcolor[rgb]{0.75,0.75,0.75}{0}&\textcolor[rgb]{0.75,0.75,0.75}{0}&\textcolor[rgb]{0.75,0.75,0.75}{0}&\textcolor[rgb]{0.75,0.75,0.75}{0}&\textcolor[rgb]{0.75,0.75,0.75}{0}&\textcolor[rgb]{0.75,0.75,0.75}{0}&$   \mu_{6+1}^\Alt \hspace{-2pt} $&\textcolor[rgb]{0.75,0.75,0.75}{0}\\[1.5pt] \hline
			  \diagbox[dir=SW,innerwidth=0.4cm]{$r$}{$q$}      &$0$  &$1$&$2$&$3$&$4$&$5$&$6$&$7$
		\end{tabular}}
      \caption{$(n,p)=(5,6)$}
      \end{subfigure}
		\begin{subfigure}{0.3\textwidth}
\scalebox{0.81}{
\begin{tabular}{c | c@{\hspace{1.3\tabcolsep}}  c@{\hspace{1.3\tabcolsep}}  c@{\hspace{1.3\tabcolsep}}   c@{\hspace{1.3\tabcolsep}}   c@{\hspace{1.3\tabcolsep}}   c@{\hspace{1.3\tabcolsep}}}
		$6$&\textcolor[rgb]{0.75,0.75,0.75}{0}&$		\mu_2^\Alt\hspace{-2pt}	 $&\textcolor[rgb]{0.75,0.75,0.75}{0}&\textcolor[rgb]{0.75,0.75,0.75}{0}&\textcolor[rgb]{0.75,0.75,0.75}{0}&\textcolor[rgb]{0.75,0.75,0.75}{0}\\[1.5pt]
		$5$&\textcolor[rgb]{0.75,0.75,0.75}{0}&\textcolor[rgb]{0.75,0.75,0.75}{0}&\textcolor[rgb]{0.75,0.75,0.75}{0}&\textcolor[rgb]{0.75,0.75,0.75}{0}&\textcolor[rgb]{0.75,0.75,0.75}{0}&\textcolor[rgb]{0.75,0.75,0.75}{0}\\[1.5pt]
		$4$&\textcolor[rgb]{0.75,0.75,0.75}{0}&\textcolor[rgb]{0.75,0.75,0.75}{0}&$   \mu_3^\Alt\hspace{-2pt}    $&\textcolor[rgb]{0.75,0.75,0.75}{0}&\textcolor[rgb]{0.75,0.75,0.75}{0}&\textcolor[rgb]{0.75,0.75,0.75}{0}\\[1.5pt]
		$3$&\textcolor[rgb]{0.75,0.75,0.75}{0}&\textcolor[rgb]{0.75,0.75,0.75}{0}&\textcolor[rgb]{0.75,0.75,0.75}{0}&\textcolor[rgb]{0.75,0.75,0.75}{0}&\textcolor[rgb]{0.75,0.75,0.75}{0}&\textcolor[rgb]{0.75,0.75,0.75}{0}\\[1.5pt]
		$2$&\textcolor[rgb]{0.75,0.75,0.75}{0}&\textcolor[rgb]{0.75,0.75,0.75}{0}&\textcolor[rgb]{0.75,0.75,0.75}{0}&$    \mu_4^\Alt\hspace{-2pt}	   $&\textcolor[rgb]{0.75,0.75,0.75}{0}&\textcolor[rgb]{0.75,0.75,0.75}{0}\\[1.5pt]
		$1$&\textcolor[rgb]{0.75,0.75,0.75}{0}&\textcolor[rgb]{0.75,0.75,0.75}{0}&\textcolor[rgb]{0.75,0.75,0.75}{0}&\textcolor[rgb]{0.75,0.75,0.75}{0}&\textcolor[rgb]{0.75,0.75,0.75}{0}&\textcolor[rgb]{0.75,0.75,0.75}{0}\\[1.5pt] 
		$0$&\textcolor[rgb]{0.75,0.75,0.75}{0}&\textcolor[rgb]{0.75,0.75,0.75}{0}&\textcolor[rgb]{0.75,0.75,0.75}{0}&\textcolor[rgb]{0.75,0.75,0.75}{0}&$ 	\mu_{4+1}^\Alt\hspace{-2pt}		$&\textcolor[rgb]{0.75,0.75,0.75}{0}\\[1.5pt] \hline
			  \diagbox[dir=SW,innerwidth=0.4cm]{$r$}{$q$}      &$0$  &$1$&$2$&$3$&$4$&$5$
		\end{tabular}}
      \caption{$(n,p)=(7,9)$}
      \end{subfigure}
		
		\begin{subfigure}{0.3\textwidth}
\scalebox{0.81}{
\begin{tabular}{c | c@{\hspace{1.3\tabcolsep}}  c@{\hspace{1.3\tabcolsep}}  c@{\hspace{1.3\tabcolsep}}  c@{\hspace{1.3\tabcolsep}}  c@{\hspace{1.3\tabcolsep}}  c@{\hspace{1.3\tabcolsep}}}
		$9$&\textcolor[rgb]{0.75,0.75,0.75}{0}&\textcolor[rgb]{0.75,0.75,0.75}{0}&\textcolor[rgb]{0.75,0.75,0.75}{0}&\textcolor[rgb]{0.75,0.75,0.75}{0}&\textcolor[rgb]{0.75,0.75,0.75}{0}&\textcolor[rgb]{0.75,0.75,0.75}{0}\\[1.5pt]
		$8$&\textcolor[rgb]{0.75,0.75,0.75}{0}&$		\mu_2^\Alt\hspace{-2pt}	 $&\textcolor[rgb]{0.75,0.75,0.75}{0}&\textcolor[rgb]{0.75,0.75,0.75}{0}&\textcolor[rgb]{0.75,0.75,0.75}{0}&\textcolor[rgb]{0.75,0.75,0.75}{0}\\[1.5pt]
		$7$&\textcolor[rgb]{0.75,0.75,0.75}{0}&\textcolor[rgb]{0.75,0.75,0.75}{0}&\textcolor[rgb]{0.75,0.75,0.75}{0}&\textcolor[rgb]{0.75,0.75,0.75}{0}&\textcolor[rgb]{0.75,0.75,0.75}{0}&\textcolor[rgb]{0.75,0.75,0.75}{0}\\[1.5pt]
		$6$&\textcolor[rgb]{0.75,0.75,0.75}{0}&\textcolor[rgb]{0.75,0.75,0.75}{0}&\textcolor[rgb]{0.75,0.75,0.75}{0}&\textcolor[rgb]{0.75,0.75,0.75}{0}&\textcolor[rgb]{0.75,0.75,0.75}{0}&\textcolor[rgb]{0.75,0.75,0.75}{0}\\[1.5pt]
		$5$&\textcolor[rgb]{0.75,0.75,0.75}{0}&\textcolor[rgb]{0.75,0.75,0.75}{0}&$     \mu_3^\Alt\hspace{-2pt}   $&\textcolor[rgb]{0.75,0.75,0.75}{0}&\textcolor[rgb]{0.75,0.75,0.75}{0}&\textcolor[rgb]{0.75,0.75,0.75}{0}\\[1.5pt]
		$4$&\textcolor[rgb]{0.75,0.75,0.75}{0}&\textcolor[rgb]{0.75,0.75,0.75}{0}&\textcolor[rgb]{0.75,0.75,0.75}{0}&\textcolor[rgb]{0.75,0.75,0.75}{0}&\textcolor[rgb]{0.75,0.75,0.75}{0}&\textcolor[rgb]{0.75,0.75,0.75}{0}\\[1.5pt]
		$3$&\textcolor[rgb]{0.75,0.75,0.75}{0}&\textcolor[rgb]{0.75,0.75,0.75}{0}&\textcolor[rgb]{0.75,0.75,0.75}{0}&\textcolor[rgb]{0.75,0.75,0.75}{0}&\textcolor[rgb]{0.75,0.75,0.75}{0}&\textcolor[rgb]{0.75,0.75,0.75}{0}\\[1.5pt]
		$2$&\textcolor[rgb]{0.75,0.75,0.75}{0}&\textcolor[rgb]{0.75,0.75,0.75}{0}&\textcolor[rgb]{0.75,0.75,0.75}{0}&$    \mu_4^\Alt\hspace{-2pt}    $&\textcolor[rgb]{0.75,0.75,0.75}{0}&\textcolor[rgb]{0.75,0.75,0.75}{0}\\[1.5pt]
		$1$&\textcolor[rgb]{0.75,0.75,0.75}{0}&\textcolor[rgb]{0.75,0.75,0.75}{0}&\textcolor[rgb]{0.75,0.75,0.75}{0}&\textcolor[rgb]{0.75,0.75,0.75}{0}&\textcolor[rgb]{0.75,0.75,0.75}{0}&\textcolor[rgb]{0.75,0.75,0.75}{0}\\[1.5pt] 
		$0$&\textcolor[rgb]{0.75,0.75,0.75}{0}&\textcolor[rgb]{0.75,0.75,0.75}{0}&\textcolor[rgb]{0.75,0.75,0.75}{0}&\textcolor[rgb]{0.75,0.75,0.75}{0}&$ 	\mu_{4+1}^\Alt\hspace{-2pt}		$&\textcolor[rgb]{0.75,0.75,0.75}{0}\\[1.5pt]  \hline
			  \diagbox[dir=SW,innerwidth=0.4cm]{$r$}{$q$}      &$0$  &$1$&$2$&$3$&$4$&$5$
		\end{tabular}}
      \caption{$(n,p)=(10,13)$}
      \end{subfigure}
		\begin{subfigure}{0.3\textwidth}
\scalebox{0.81}{
\begin{tabular}{c | c@{\hspace{1.3\tabcolsep}}  c@{\hspace{1.3\tabcolsep}}  c@{\hspace{1.3\tabcolsep}}  c@{\hspace{1.3\tabcolsep}}  c@{\hspace{1.3\tabcolsep}}  c@{\hspace{1.3\tabcolsep}}}
		$9$&\textcolor[rgb]{0.75,0.75,0.75}{0}&\textcolor[rgb]{0.75,0.75,0.75}{0}&\textcolor[rgb]{0.75,0.75,0.75}{0}&\textcolor[rgb]{0.75,0.75,0.75}{0}&\textcolor[rgb]{0.75,0.75,0.75}{0}&\textcolor[rgb]{0.75,0.75,0.75}{0}\\[1.5pt]
		$8$&\textcolor[rgb]{0.75,0.75,0.75}{0}&\textcolor[rgb]{0.75,0.75,0.75}{0}&\textcolor[rgb]{0.75,0.75,0.75}{0}&\textcolor[rgb]{0.75,0.75,0.75}{0}&\textcolor[rgb]{0.75,0.75,0.75}{0}&\textcolor[rgb]{0.75,0.75,0.75}{0}\\[1.5pt]
		$7$&\textcolor[rgb]{0.75,0.75,0.75}{0}&$		\mu_2^\Alt\hspace{-2pt}	 $&\textcolor[rgb]{0.75,0.75,0.75}{0}&\textcolor[rgb]{0.75,0.75,0.75}{0}&\textcolor[rgb]{0.75,0.75,0.75}{0}&\textcolor[rgb]{0.75,0.75,0.75}{0}\\[1.5pt]
		$6$&\textcolor[rgb]{0.75,0.75,0.75}{0}&\textcolor[rgb]{0.75,0.75,0.75}{0}&\textcolor[rgb]{0.75,0.75,0.75}{0}&\textcolor[rgb]{0.75,0.75,0.75}{0}&\textcolor[rgb]{0.75,0.75,0.75}{0}&\textcolor[rgb]{0.75,0.75,0.75}{0}\\[1.5pt]
		$5$&\textcolor[rgb]{0.75,0.75,0.75}{0}&\textcolor[rgb]{0.75,0.75,0.75}{0}&\textcolor[rgb]{0.75,0.75,0.75}{0}&\textcolor[rgb]{0.75,0.75,0.75}{0}&\textcolor[rgb]{0.75,0.75,0.75}{0}&\textcolor[rgb]{0.75,0.75,0.75}{0}\\[1.5pt]
		$4$&\textcolor[rgb]{0.75,0.75,0.75}{0}&\textcolor[rgb]{0.75,0.75,0.75}{0}&$   \mu_3^\Alt\hspace{-2pt}    $&\textcolor[rgb]{0.75,0.75,0.75}{0}&\textcolor[rgb]{0.75,0.75,0.75}{0}&\textcolor[rgb]{0.75,0.75,0.75}{0}\\[1.5pt]
		$3$&\textcolor[rgb]{0.75,0.75,0.75}{0}&\textcolor[rgb]{0.75,0.75,0.75}{0}&\textcolor[rgb]{0.75,0.75,0.75}{0}&\textcolor[rgb]{0.75,0.75,0.75}{0}&\textcolor[rgb]{0.75,0.75,0.75}{0}&\textcolor[rgb]{0.75,0.75,0.75}{0}\\[1.5pt]
		$2$&\textcolor[rgb]{0.75,0.75,0.75}{0}&\textcolor[rgb]{0.75,0.75,0.75}{0}&\textcolor[rgb]{0.75,0.75,0.75}{0}&\textcolor[rgb]{0.75,0.75,0.75}{0}&\textcolor[rgb]{0.75,0.75,0.75}{0}&\textcolor[rgb]{0.75,0.75,0.75}{0}\\[1.5pt]
		$1$&\textcolor[rgb]{0.75,0.75,0.75}{0}&\textcolor[rgb]{0.75,0.75,0.75}{0}&\textcolor[rgb]{0.75,0.75,0.75}{0}&$    \mu_4^\Alt\hspace{-2pt}    $&\textcolor[rgb]{0.75,0.75,0.75}{0}&\textcolor[rgb]{0.75,0.75,0.75}{0}\\[1.5pt] 
		$0$&\textcolor[rgb]{0.75,0.75,0.75}{0}&\textcolor[rgb]{0.75,0.75,0.75}{0}&\textcolor[rgb]{0.75,0.75,0.75}{0}&\textcolor[rgb]{0.75,0.75,0.75}{0}&$ 	\mu_{4+1}^\Alt\hspace{-2pt}		$&\textcolor[rgb]{0.75,0.75,0.75}{0}\\[1.5pt]  \hline
			  \diagbox[dir=SW,innerwidth=0.4cm]{$r$}{$q$}      &$0$  &$1$&$2$&$3$&$4$&$5$
		\end{tabular}}
      \caption{$(n,p)=(9,12)$}
      \end{subfigure}
\caption{Possibly non-zero entries, with its rank, of the $E^\infty$-page of the spectral sequence $E^1_{r,q}= H^{\Alt_{r+1}}_q\left(D^{r+1}\left(F\right),D^{r+1}\left(f_s\right)\right)$ for a map germ $f:\GS{n}{p}$. Notice the shift $H_{m+1}(\im F ,\im f_s)\cong H_m(\im f_s)$ when $m>0$.}
\label{fig: las mejores tablas del universo}
\end{figure}

\begin{remark}\label{rem: expected dim}
For non-degenerated germs, the dimensions where $\im f_s$ has (possibly) non-trivial homology are (cf. \cref{fig: las mejores tablas del universo})
\begin{itemize}
	\item $d_k+k-1$ for $d_k\geq0$, from $\beta^\Alt\big(D^k(f_s)\big)$ or, equivalently, from $\mu_k^\Alt(f)$; and
   \item $\kappa=\left\lfloor \frac{p}{p-n}\right\rfloor-1$, from $\mu_{\kappa+1}^\Alt(f)$. 
\end{itemize}
\end{remark}

\begin{figure}[htb]
\subfloat[$(n,p)=(2,4)$]{
\begin{tabular}{c | c@{\hspace{1.3\tabcolsep}}  c@{\hspace{1.3\tabcolsep}}  c@{\hspace{1.3\tabcolsep}}   c@{\hspace{1.3\tabcolsep}}   c@{\hspace{1.3\tabcolsep}}   c@{\hspace{1.3\tabcolsep}}}
		$4	 $&\textcolor[rgb]{0.75,0.75,0.75}{0}&\textcolor[rgb]{0.75,0.75,0.75}{0}&\textcolor[rgb]{0.75,0.75,0.75}{0}&\textcolor[rgb]{0.75,0.75,0.75}{0}&\textcolor[rgb]{0.75,0.75,0.75}{0}&\textcolor[rgb]{0.75,0.75,0.75}{0}\\[1.5pt]
		$3	 $&\textcolor[rgb]{0.75,0.75,0.75}{0}&\textcolor[rgb]{0.75,0.75,0.75}{0}&\textcolor[rgb]{0.75,0.75,0.75}{0}&\textcolor[rgb]{0.75,0.75,0.75}{0}&\textcolor[rgb]{0.75,0.75,0.75}{0}&\textcolor[rgb]{0.75,0.75,0.75}{0}\\[1.5pt]
		$2	 $&\textcolor[rgb]{0.75,0.75,0.75}{0}&\textcolor[rgb]{0.75,0.75,0.75}{0}&\textcolor[rgb]{0.75,0.75,0.75}{0}&\textcolor[rgb]{0.75,0.75,0.75}{0}&\textcolor[rgb]{0.75,0.75,0.75}{0}&\textcolor[rgb]{0.75,0.75,0.75}{0}\\[1.5pt]
		$1	 $&\textcolor[rgb]{0.75,0.75,0.75}{0}&$   \mu_{2}^\Alt \hspace{-2pt} $&\textcolor[rgb]{0.75,0.75,0.75}{0}&\textcolor[rgb]{0.75,0.75,0.75}{0}&\textcolor[rgb]{0.75,0.75,0.75}{0}&\textcolor[rgb]{0.75,0.75,0.75}{0}\\[1.5pt] 
		$0	 $&\textcolor[rgb]{0.75,0.75,0.75}{0}&\textcolor[rgb]{0.75,0.75,0.75}{0}&$   \mu_{2+1}^\Alt \hspace{-2pt} $&\textcolor[rgb]{0.75,0.75,0.75}{0}&\textcolor[rgb]{0.75,0.75,0.75}{0}&\textcolor[rgb]{0.75,0.75,0.75}{0}\\[1.5pt] \hline
			  \diagbox[dir=SW,innerwidth=0.4cm]{$r$}{$q$}      &$0$  &$1$&$2$&$3$&$4$&$5$
		\end{tabular}}
		\hspace{1cm}
		\subfloat[$(n,p)=(2,5)$]{
\begin{tabular}{c | c@{\hspace{1.3\tabcolsep}}  c@{\hspace{1.3\tabcolsep}}  c@{\hspace{1.3\tabcolsep}}   c@{\hspace{1.3\tabcolsep}}   c@{\hspace{1.3\tabcolsep}}   c@{\hspace{1.3\tabcolsep}}   }
		$4	 $&\textcolor[rgb]{0.75,0.75,0.75}{0}&\textcolor[rgb]{0.75,0.75,0.75}{0}&\textcolor[rgb]{0.75,0.75,0.75}{0}&\textcolor[rgb]{0.75,0.75,0.75}{0}&\textcolor[rgb]{0.75,0.75,0.75}{0}&\textcolor[rgb]{0.75,0.75,0.75}{0}\\[1.5pt]
		$3	 $&\textcolor[rgb]{0.75,0.75,0.75}{0}&\textcolor[rgb]{0.75,0.75,0.75}{0}&\textcolor[rgb]{0.75,0.75,0.75}{0}&\textcolor[rgb]{0.75,0.75,0.75}{0}&\textcolor[rgb]{0.75,0.75,0.75}{0}&\textcolor[rgb]{0.75,0.75,0.75}{0}\\[1.5pt]
		$2	 $&\textcolor[rgb]{0.75,0.75,0.75}{0}&\textcolor[rgb]{0.75,0.75,0.75}{0}&\textcolor[rgb]{0.75,0.75,0.75}{0}&\textcolor[rgb]{0.75,0.75,0.75}{0}&\textcolor[rgb]{0.75,0.75,0.75}{0}&\textcolor[rgb]{0.75,0.75,0.75}{0}\\[1.5pt]
		$1	 $&\textcolor[rgb]{0.75,0.75,0.75}{0}&\textcolor[rgb]{0.75,0.75,0.75}{0}&\textcolor[rgb]{0.75,0.75,0.75}{0}&\textcolor[rgb]{0.75,0.75,0.75}{0}&\textcolor[rgb]{0.75,0.75,0.75}{0}&\textcolor[rgb]{0.75,0.75,0.75}{0}\\[1.5pt] 
		$0	 $&\textcolor[rgb]{0.75,0.75,0.75}{0}&$   \mu_{1+1}^\Alt \hspace{-2pt} $&\textcolor[rgb]{0.75,0.75,0.75}{0}&\textcolor[rgb]{0.75,0.75,0.75}{0}&\textcolor[rgb]{0.75,0.75,0.75}{0}&\textcolor[rgb]{0.75,0.75,0.75}{0}\\[1.5pt] \hline
			  \diagbox[dir=SW,innerwidth=0.4cm]{$r$}{$q$}      &$0$  &$1$&$2$&$3$&$4$&$5$
		\end{tabular}}
\caption{Comparison between the spectral sequences of non-degenerated and degenerated map germs $f:\GS{n}{p}$ (where the stable perturbation only has 0-homology, i.e., the expected dimension $d_2$ is negative). Notice the equality $\beta_{1}\big(\im(F),\im(f_s)\big)+1= \beta_0\big(\im(f_s)\big)$ in the degenerated case.}
\label{fig: las segundas mejores tablas del universo}
\end{figure}

Now we can give a workable expression for the image Milnor number $\mu_I$.

\begin{theorem}\label{thm: conustedesmui}
Given an unstable $\eqA$-finite germ $f:\GS{n}{p}$ of corank one, $n<p\leq2n$,
\begin{equation*}
	\mu_I(f)=\sum_{k=2}^{d(f)}  \frac{1}{k!}\left(\sum_{ d_k^\sigma\geq0}\mu\big(D^k(f)^\sigma\big)-\sum_{ d_k^\sigma<0}(-1)^{d_k^\sigma}\beta_0\big(D^k(f)^\sigma\big)\right)+	\genfrac(){0pt}{0}{s(f)-1}{d(f)},
	\end{equation*}
where 
$ 	\genfrac(){0pt}{0}{s(f)-1}{d(f)}=0 $
if $s(f)\leq d(f)$.
\end{theorem}
\begin{proof}
The proof follows from \cref{conjrgc,lem: suma mukalt}.  
\end{proof}

\begin{remark}\label{rem: marar-mond formulae}
One can give a similar version for $\nu_I(f)$. Indeed, \cref{thm: conustedesmui} (together with \cref{eqthm: mu+-} of \cref{conjrgc}) should be regarded as the generalization of the Marar-Mond formulae (see \cite{Marar1989}, in particular Propositions 3.2 and 3.3 and Theorem 3.4, cf. \cite{Mond1991}), which relate the image Milnor number $\mu_I$ of a mono-germ $f:(\CC^2,0)\to(\CC^3,0)$ with the number of triple points $T$, cross-caps $C$ (in $\im f_s$) and the Milnor number of the curve of double points $D(f)\subset \CC^2$. A generalization of these formulae was given, only in terms of Euler characteristics (and for mono-germs), by Marar in \cite[Theorem 3.1]{Marar1991}, using a combinatorial argument. This argument, however, leaves inaccessible $\mu_I$ because it relies on Euler-Poincaré characteristics. The version of \cref{thm: conustedesmui} for $\nu_I$ would be a reformulation of this generalization by Marar, with a different proof.
\end{remark}

The ideas that led to \cref{thm: conustedesmui} allow us to prove a modification of the weak version of the Mond conjecture for corank one and $n<p$. This is proved for corank one in \cite[Theorem 3.9]{GimenezConejero2022} and for any corank in \cite[Theorem 2.15]{GimenezConejero2023a} when $p=n+1$.

\begin{theorem}\label{wmcnp}
Let $f\colon\CCS{n}\rightarrow \CCzero{p}$ be $\eqA$-finite of corank 1, $n<p$. Then, $\mu_I(f)=0$ if, and only if, $f$ is stable or, alternatively, strongly contractible.
\end{theorem}
\begin{proof}
This follows from \cref{conjrgc,cor: sing iff alt,lem: suma mukalt}.  
\end{proof}

\begin{remark}
The germ given in \cref{ex:disgusting} illustrates the need of adding the strongly contractible case. Strongly contractible germs are, precisely, counterexamples to the original formulation of the result: $\mu_I(f)=0$ if, and only if, $f$ is stable. 

Moreover, it is easy to see why this would fail for $\nu_I(f)$, as there could be some cancellation. For example, the germ $f:(\CC^4,0)\to(\CC^6,0)$ such that
$$ f(x_1,x_2,x_3,y)=(x_1,x_2,x_3,\ y^3+x_1^2y,\ y^4+x_2y,\ y^5+x_3y) $$
was given in \cite[Example 4.26]{Houston1997} and the $\im f_s$ ($f_s$ stable) has the homotopy type of two $2$-spheres and two $3$-spheres, so $\mu_I(f)=4$ but $\nu_I(f)=0$ being $f$ unstable (and not strongly contractible).
\end{remark}

\subsection{Unexpected homology}\label{subs:unexhom}

The previous definitions and results deal with the homology of images of stable perturbations, but one can wonder what happens with a deformation $f_t$ that is not necessarily locally stable. In this section we develop what we call \textsl{unexpected homology}. It is homology of $\im(f_t)$ that is killed when we deform further $f_t$ to a locally stable map. 
\newline

 For $p=n+1$ we have a conservation principle proved in \cite[Theorem 2.8]{GimenezConejero2022} using a result for hypersurfaces (in particular, \cite[Theorem 2.3]{Siersma1991} together with \cite[Lemma 2.4]{GimenezConejero2022}). As this result is only for hypersurfaces, the same proof cannot be used when $p>n+1$. 

\begin{figure}
	\centering
		\includegraphics[width=0.70\textwidth]{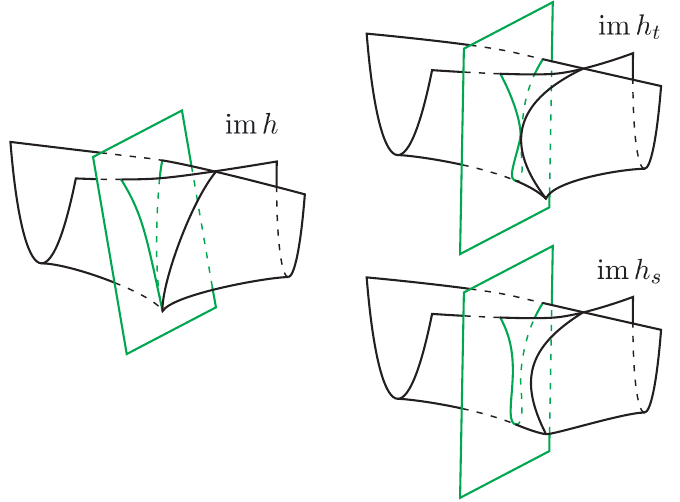}
	\caption{Creation of unexpected homology in $h_t$ and its annihilation in a locally stable perturbation $h_s$.}
	\label{fig:UnexpectedHomology}
\end{figure}

Indeed, in light of \cref{ex:multiSS}, the conservation principle of $\mu_I$ or $\nu_I$ is not true in general: that deformation $h_t$ of the bi-germ $h$ has homology that does not come from any of the spaces $D^k(h_t)$ with $d_k\geq0$ or combination of branches (because $d(h)\geq s(h)$), hence 
\begin{align*}
	\mu_I(h)&<    \sum_{i>0} \beta_i\big(\im(h_t)\big)                               + \sum_{y\in\im(h_t)} \mu_I(h_t;y)          ,\text{ and }\\
	\nu_I(h)&\neq     \sum_{i>0} (-1)^{i+d_2+1}\beta_i\big(\im(h_t)\big)     + \sum_{y\in\im(h_t)} \nu_I(h_t;y).
\end{align*}
Obviously, if we keep deforming this map $h_t$ to have something locally stable, $h_s$, this \textsl{unexpected homology} is killed (see \cref{fig:UnexpectedHomology}). This follows from the \textsc{icss}, or \cref{lem: suma mukalt}. 

Let us see the \textsc{icss} for $(\im H, \im h_t)$, with $H$ a versal unfolding of $h$. The first page of the spectral sequence $E^1_{r,q}= H^{\Alt_{r+1}}_q\left(D^{r+1}\left(H\right),D^{r+1}\left(h_t\right)\right)$ is 
$$
\begin{tabular}{c | c@{\hspace{1.3\tabcolsep}}  c@{\hspace{1.3\tabcolsep}}  c@{\hspace{1.3\tabcolsep}}   c@{\hspace{1.3\tabcolsep}}   c@{\hspace{1.3\tabcolsep}}   c@{\hspace{1.3\tabcolsep}}}
		$H^\Alt_3$&\textcolor[rgb]{0.75,0.75,0.75}{0}&\textcolor[rgb]{0.75,0.75,0.75}{0}   &\textcolor[rgb]{0.75,0.75,0.75}{0}&\textcolor[rgb]{0.75,0.75,0.75}{0}\\[2pt]
		$H^\Alt_2$&\textcolor[rgb]{0.75,0.75,0.75}{0}&$\mathop{\scalebox{1.5}{$\ast$}}$&\textcolor[rgb]{0.75,0.75,0.75}{0}&\textcolor[rgb]{0.75,0.75,0.75}{0}\\[2pt]
		$H^\Alt_1$&\textcolor[rgb]{0.75,0.75,0.75}{0}&\textcolor[rgb]{0.75,0.75,0.75}{0}   &$\ZZ$                                            &\textcolor[rgb]{0.75,0.75,0.75}{0}\\[2pt] 
		$H^\Alt_0$&\textcolor[rgb]{0.75,0.75,0.75}{0}&\textcolor[rgb]{0.75,0.75,0.75}{0}   &\textcolor[rgb]{0.75,0.75,0.75}{0}&\textcolor[rgb]{0.75,0.75,0.75}{0}\\[2pt] \hline
			  \diagbox[dir=SW,innerwidth=0.8cm]{$q$}{$r$}      &$D^{0+1}$  &$D^{1+1}$&$D^{2+1}$&$D^{3+1}$
		\end{tabular}$$
We can compare it to the usual \textsc{icss} (after collapsing) for $(\im A, \im a_s)$ of an arbitrary multi-germ $a$ in the same dimensions with possibly $s(a)>d(a)$, i.e., the $E^\infty$-page of the spectral sequence $E^1_{r,q}= H^{\Alt_{r+1}}_q\left(D^{r+1}\left(A\right),D^{r+1}\left(a_s\right)\right)$: 
$$
\begin{tabular}{c | c@{\hspace{1.3\tabcolsep}}  c@{\hspace{1.3\tabcolsep}}  c@{\hspace{1.3\tabcolsep}}   c@{\hspace{1.3\tabcolsep}}   c@{\hspace{1.3\tabcolsep}}   c@{\hspace{1.3\tabcolsep}}}
		$H^\Alt_3$&\textcolor[rgb]{0.75,0.75,0.75}{0}&\textcolor[rgb]{0.75,0.75,0.75}{0}   &\textcolor[rgb]{0.75,0.75,0.75}{0}&\textcolor[rgb]{0.75,0.75,0.75}{0}\\[2pt]
		$H^\Alt_2$&\textcolor[rgb]{0.75,0.75,0.75}{0}&$\mathop{\scalebox{1.5}{$\ast$}}$&\textcolor[rgb]{0.75,0.75,0.75}{0}&\textcolor[rgb]{0.75,0.75,0.75}{0}\\[2pt]
		$H^\Alt_1$&\textcolor[rgb]{0.75,0.75,0.75}{0}&\textcolor[rgb]{0.75,0.75,0.75}{0}   &\textcolor[rgb]{0.75,0.75,0.75}{0} &\textcolor[rgb]{0.75,0.75,0.75}{0}\\[2pt] 
		$H^\Alt_0$&\textcolor[rgb]{0.75,0.75,0.75}{0}&\textcolor[rgb]{0.75,0.75,0.75}{0}   &$\mathop{\scalebox{1.5}{$\ast$}}$&\textcolor[rgb]{0.75,0.75,0.75}{0}\\[2pt] \hline
			  \diagbox[dir=SW,innerwidth=0.8cm]{$q$}{$r$}      &$D^{0+1}$  &$D^{1+1}$&$D^{2+1}$&$D^{3+1}$
		\end{tabular}$$
		
This \textsl{unexpected homology} of \cref{ex:multiSS} comes from alternating homology of $D^3(h_t)$ that cannot be found on $D^3(H)$, which is equivalent to having non-trivial $$H_1^\Alt\big(D^3(H),D^3(h_t)\big), $$
by the exact sequence of the pair in alternating homology. 
\newline

Let us explain this in detail. 

For the moment, let us assume that we are in some dimensions where there is no $k$ such that $d_k=0$ (as it is the case of $h$ in \cref{ex:multiSS}). 
 In general, for arbitrary $f$, we could have several spaces $D^k(f_t)$ with $d_k<0$ such that they have alternating homology that is not found in $D^k(F)$, or, equivalently, non-trivial $H_1^\Alt$ of the pair. These homologies 
$$ H_1^\Alt\big(D^k(F),D^k(f_t)\big) $$
 appear in the \textsc{icss} in the first row of the first page, $E^1_{*,1}(F, f_t)$, and it collapses in the second page leaving a possibly non-trivial term in the column where $D^{\kappa+1}$ goes, i.e., on the entry $E^2_{\kappa,1}(F,f_t)$, where $\kappa=\left\lfloor \frac{p}{p-n}\right\rfloor$. For example, if $\ell$ is the maximal multiplicity within the spaces $D^{k}(f_t)$ that have
$$ H_1^\Alt\big(D^k(F),D^k(f_t)\big)\ncong 0,$$
we have an exact sequence within the first row of the first page of the 
 \textsc{icss}  with the following combinatorial shape:
\begin{equation}\label{eq: sequence unexpected}
\ZZ^{\genfrac(){0pt}{1}{\ell}{\kappa+1}}\leftarrow \ZZ^{\genfrac(){0pt}{1}{\ell}{\kappa+2}} \leftarrow\cdots\leftarrow\ZZ^{\genfrac(){0pt}{1}{\ell}{\ell-1}}\leftarrow\ZZ^{\genfrac(){0pt}{1}{\ell}{\ell}}\leftarrow0,\end{equation}
followed by zeroes. Of course, we have more of these sequences as direct sums on the same entries of the \textsc{icss}, depending on the orbits and projections of the $D^k(f_t)$ with $d_k<0$. Indeed, this is the same combinatorial reasoning that is used to prove that the $0$-th row of the first page of the \textsc{icss} $E^r_{*,*}(F,f_s)$ forms an exact sequence and all the homology collapses to give the unique (possibly) non-trivial $E^2_{\kappa+1,0}$, which defines $\mu_{d(f)+1}^\Alt$. This is proven in \cite[Lemma 3.3]{Houston2010}  for rational homology, which is equivalent to integer homology in this case because we only deal with some $H_0$ of multiple points and the combinatorics of the projections between them (see also \cite[Figure 2]{Houston2010} for a better insight).  
So we see that the term is non-trivial if there is such alternating homology for some $k$. See \cite[Lemmas 3.3 and 3.4]{Houston2010} for more details. 

In contrast, if $n$ and $p$ are such that $d_\kappa=0$, all the sequences of the shape of \cref{eq: sequence unexpected} have an extra term coming from $D^\kappa(f_t)$:
\begin{equation}\label{eq: sequence no unexpected}
 \ZZ^{\genfrac(){0pt}{1}{\ell}{\kappa}}\leftarrow\ZZ^{\genfrac(){0pt}{1}{\ell}{\kappa+1}}\leftarrow \ZZ^{\genfrac(){0pt}{1}{\ell}{\kappa+2}} \leftarrow\cdots\leftarrow\ZZ^{\genfrac(){0pt}{1}{\ell}{\ell-1}}\leftarrow\ZZ^{\genfrac(){0pt}{1}{\ell}{\ell}}\leftarrow0.
\end{equation}
With 
 the same argument, we see that this collapses to give a (possibly) non-trivial $E^2_{\kappa-1,1}(F,f_t)$, in the column of $D^\kappa$, and always killing $E^2_{\kappa,1}(F,f_t)$ (contributing trivially to the homology of $\im(f_t)$ in the corresponding rank). 

\begin{definition}\label{def: unexpected homology}
We say that a perturbation $f_t$ of an $\eqA$-finite map germ $f:(\CC^n,S)\to(\CC^p,0)$, $n<p$, has \textit{unexpected homology} if 
there are $k$ such that 
$$ H_1^\Alt\big( D^k(F), D^k(f_t) \big)\ncong0, $$
where $F$ is a stable unfolding, and all satisfy that $d_k(f)<0$ (in contrast with the situation of \cref{eq: sequence no unexpected}). In that case, the \textit{unexpected homology} is found in a unique rank:
$$ H_{\kappa} (\im f_t) \ncong 0, $$
with $\kappa=\left\lfloor \frac{p}{p-n}\right\rfloor$.
\end{definition}

\begin{remark}\label{rem: unexpected}
An $\eqA$-finite map germ $f:(\CC^n,S)\to(\CC^p,0)$, $n<p$, such that $D^k(f)=\varnothing$ for every negative expected dimension $d_k$ or such that $\frac{p}{p-n}\in\ZZ^+$, has no deformations with unexpected homology. 

To understand \textsl{what unexpected homology is}, and not only \textsl{where it is}, recall the previous discussion and \cref{rem: expected dim}: if there are no $k$ such that $d_k+k-1=\kappa$, then the unexpected homology of $\im(f_t)$ is $$ H_{\kappa} (\im f_t),$$
if this is not trivial. 

Indeed, the only case where we have a $k$ such that $d_k+k-1=\kappa$ is when $d_\kappa= 1$, in which case such $k$ is necessarily $\kappa$. This is because $d_k+k-1>d_\kappa+\kappa-1$ and $d_k+k-1>k$ if, simultaneously, $p>n+1$ and $k<\kappa$ (so $d_k>0$). In that case, the unexpected homology is given by the difference 
$$ \beta_\kappa(\im f_t) - \beta_1^\Alt \big( D^\kappa(f_t) \big),$$
if it is positive.
\end{remark}

Now, we have the best approximation to a conservation principle for $\mu_I$ and $\nu_I$: we have to eliminate the unexpected homology when necessary. This generalises  \cite[Theorem 2.8]{GimenezConejero2022}.

\begin{theorem}\label{thm: conservation muinui}
Let $f:\GS{n}{p}$ be an $\eqA$-finite germ  of corank one, $n<p\leq2n$. Then, for any perturbation $f_t$, if $\frac{p}{p-n}\notin\ZZ^+$:
\begin{equation}\label{eq: conservation exceptional}
\begin{aligned}
	\mu_I(f)&=    \sum_{i\neq0,\kappa} \beta_{i}(\im f_t)                               + \sum_{y\in\im(f_t)} \mu_I(f_t;y)         - \beta_\kappa(\im f_t) +\delta ,\text{ and }\\
	\nu_I(f)&=     \sum_{i\neq0,\kappa} (-1)^{i+d_2+1}\beta_{i}(\im f_t)     +             \sum_{y\in\im(f_t)} \nu_I(f_t;y) - (-1)^{\kappa+d_2+1}\big(\beta_\kappa(\im f_t)-\delta\big);
\end{aligned}
\end{equation}
where $\delta=\beta_1^\Alt\big(D^\kappa(f_t)\big)$ if $d_\kappa=1$ and $0$ otherwise (see \cref{rem: unexpected}).

If $\frac{p}{p-n}\in\ZZ^+$:
\begin{equation}\label{eq: conservation no exceptional}
\begin{aligned}
	\mu_I(f)&=    \sum_{i\neq0} \beta_{i}(\im f_t)                               + \sum_{y\in\im(f_t)} \mu_I(f_t;y)   ,\text{ and }\\
	\nu_I(f)&=     \sum_{i\neq0} (-1)^{i+d_2+1}\beta_{i}(\im f_t)     +             \sum_{y\in\im(f_t)} \nu_I(f_t;y),
\end{aligned}
\end{equation}
\end{theorem}
\begin{proof}
Assume that $\frac{p}{p-n}\notin\ZZ^+$ and let us focus on the case of $\mu_I$, the case $\nu_I$ is analogous keeping track of the signs.

Recall that (\cref{rem: expected dim}), in general, $\mu_k^\Alt$ contributes in rank $d_k+k-1$ for $k\leq\kappa$ and $\mu_{\kappa+1}^\Alt$ in rank $\kappa-1$:
\begin{equation}\label{eq:acaba} \begin{aligned}
\mu_k^\Alt(f)  &= \beta_{d_k+k-1} (\im f_s) \quad\quad\quad\quad\text{ for } k=2,\dots,\kappa-1\\
\mu_\kappa^\Alt(f)  &=  \begin{cases}
\beta_{d_\kappa+\kappa-1} (\im f_s) & \text{ if } d_\kappa\neq0\\
\beta_{d_\kappa+\kappa-1} (\im f_s) -\mu_{\kappa+1}^\Alt(f) & \text{ if } d_\kappa=0
\end{cases}\\
\mu_{\kappa+1}^\Alt(f)  &=  \begin{cases}
\beta_{\kappa-1} (\im f_s) & \quad\quad\text{if } d_\kappa\neq0\\
\beta_{\kappa-1} (\im f_s) -\mu_{\kappa}^\Alt(f) & \quad\quad\text{if } d_\kappa=0.
\end{cases}
\end{aligned}
\end{equation}

By the conservation principle, \cref{thm: conservation mu tau f}, and \cref{lem: suma mukalt}; $\beta_k^\Alt\big(D^k(f_t)\big)$ 
 for $k\leq \kappa$ contributes positively to \cref{eq: conservation exceptional}. This explains most of the positive terms of \cref{eq: conservation exceptional}. It
 shows that
\begin{equation}
\sum_{k< \kappa}\mu_k^\Alt(f)=    \sum_{k<\kappa} \beta_{d_k+k-1}(\im f_t)                          + \sum_{y\in\im(f_t)} \sum_{k< \kappa}\mu_k^\Alt(f_t;y),
\label{eq:confusing1}
\end{equation} 
and
\begin{equation}
\mu_\kappa^\Alt(f)=
\beta_{d_\kappa}^\Alt\big(D^\kappa(f_t)\big)                    
+ \sum_{y\in\im(f_t)} \mu_\kappa^\Alt(f_t;y).
\label{eq:confusing2}
\end{equation} 
If $d_\kappa\neq0$, we also have
\[\mu_\kappa^\Alt(f)= \beta_{d_\kappa+\kappa-1}(\im f_t)+ \sum_{y\in\im(f_t)} \mu_\kappa^\Alt(f_t;y).\]
However, if $d_\kappa=0$, it can happen that $\beta_{d_\kappa+\kappa-1}(\im f_t)$ is equal to $\beta_{d_\kappa}^\Alt\big(D^\kappa(f_t)\big)   $ plus a contribution from the spaces $H_0^\Alt\big(D^k(f_t)\big)$ for $k>\kappa$  (see \cref{eq:acaba}), similarly to the case for the stable perturbation but with $\mu_\kappa^\Alt(f)$ and $\mu_{\kappa+1}^\Alt(f)$ (both contribute in rank $\kappa-1$).

Hence, to prove \cref{eq: conservation exceptional}, it is sufficient to show that
\begin{equation}
\mu_{\kappa+1}^\Alt(f)=
\begin{cases}\beta_{\kappa-1}(\im f_t) +\sum_{y\in\im(f_t)} \mu_{\kappa+1}^\Alt(f_t;y) -\beta_{\kappa}\big(\im(f_t)\big)+\delta & \text{ if }d_\kappa\neq0\\
\beta_{\kappa-1}(\im f_t) -\beta_{d_\kappa}^\Alt\big(D^\kappa(f_t)\big)+\sum_{y\in\im(f_t)} \mu_{\kappa+1}^\Alt(f_t;y) -\beta_{\kappa}\big(\im(f_t)\big)+\delta & \text{ if }d_\kappa=0
\end{cases}
\label{eq:confusing3}
\end{equation}
If we show this equation, the result follows by taking the sum of \cref{eq:confusing1,eq:confusing2,eq:confusing3}.
\newline

Finally, we use the spectral sequence to compute the terms of \cref{eq:confusing3}:
\begin{itemize}
	\item To compute $\mu_{\kappa+1}^\Alt(f)$, we need the homologies 
\[H_0^\Alt\big(D^k(f)\big)\cong H_0^\Alt\big(D^k(F)\big)\cong H_0^\Alt\big(D^k(F),D^k(f_s)\big),\]
 for a versal unfolding $F$, a stable perturbation $f_s$ and $k>\kappa$ (therefore $D^k(f_s)=\varnothing$). 

\item To compute $\beta_{\kappa-1}(\im f_t)$, or $\beta_{\kappa-1}(\im f_t) -\beta_{d_\kappa}^\Alt\big(D^\kappa(f_t)\big)$ if $d_\kappa=0$, we need the cokernels
$$ C\coloneqq \coker\left( H_0^\Alt\big(D^k(f_t)\big)\to H_0^\Alt\big(D^k(F)\big)\right)\cong H_0^\Alt\big( D^k(F), D^k(f_t) \big).$$

\item To compute $\sum_y\mu_{\kappa+1}^\Alt(f_t;y)$, we need the spaces
$$H_0^\Alt\big(D^k(f_t)\big)\cong I\oplus K$$
 with $k>\kappa$ and
\begin{align*}
I \coloneqq \im&\left( H_0^\Alt\big(D^k(f_t)\big)\to H_0^\Alt\big(D^k(F)\big)\right),\\
K\coloneqq \ker&\left( H_0^\Alt\big(D^k(f_t)\big)\to H_0^\Alt\big(D^k(F)\big)\right)\cong H_1^\Alt\big( D^k(F), D^k(f_t) \big),
\end{align*}
\item To compute $\beta_{\kappa}\big(\im(f_t)\big)+\delta$,
 which is precisely the unexpected homology, we need the spaces $K$.
\end{itemize}


Hence, since
\begin{equation}
\label{eq:CIK} \dim H_0^\Alt\big(D^k(F)\big) = \dim C + \dim I = \dim C + \dim H_0^\Alt\big(D^k(f_t)\big) - \dim K,
\end{equation}
we can show that \cref{eq:confusing3} holds, after a combinatorial argument.
\newline

The \cref{eq: conservation no exceptional} follows from the previous argument, since there is no unexpected homology.
\end{proof}

\begin{figure}
	\centering
		\includegraphics[width=0.70\textwidth]{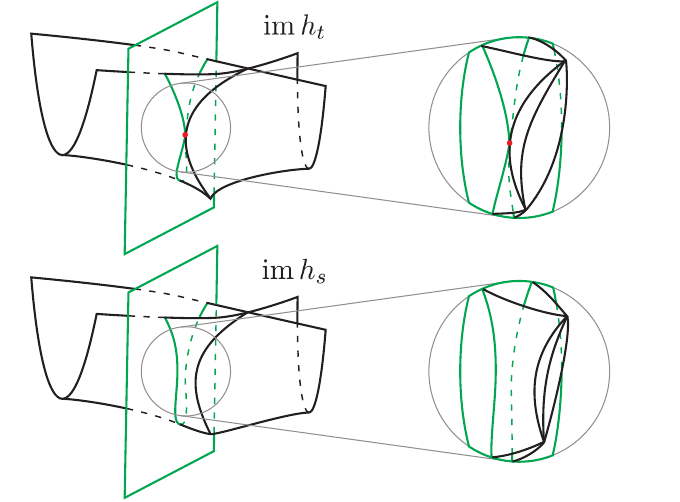}
	\caption{Comparison between unexpected homology and $\mu_{d(h)+1}^\Alt(h_t)$: the homology of $\im h_t$ is killed when we take the stable perturbation $h_s$, and this unexpected homology is reflected in the $\mu_{d(h)+1}^\Alt(h_t)$ at the new instability (red dot). See \cref{ex: killing}.}
	\label{fig:UnexpectedHomologyDetail}
\end{figure}

\begin{example}\label{ex: killing}
In the deformation $h_t$ of \cref{ex:multiSS}, there is an instability at the origin with $\mu_I(h_t;0)=1$ given only by $\mu_{d(h)+1}^\Alt(h_t)$, but also unexpected homology $\beta_{\kappa}(\im h_t)=1$. Cancelling these quantities we recover the conservation of the $\mu_I$ (this is illustrated in \cref{fig:UnexpectedHomologyDetail}).
\end{example}

\section{Houston's conjecture on excellent unfoldings}\label{sec: houston c}

Houston conjectured in \cite[Conjecture 6.2]{Houston2010} that \textsl{excellency }in families can be controlled by means of the image Milnor number, for any $p>n$, in corank one. In the same article, Houston gave a partial answer in \cite[Corollary 5.9]{Houston2010}, where the family has to be \textsl{good} (and more assumptions). 
 This conjecture was proved by the author and Nuño-Ballesteros for $p=n+1$ in \cite[Theorem 4.3]{GimenezConejero2022}. We present now a complete answer using our new approach, giving counterexamples and providing a correct modification of the statement of the conjecture which we also prove. 
 An answer to the conjecture \cite[Conjecture 4.6]{GimenezConejero2022}, which asks about the converse, is also given.

\begin{definition}\label{def: excellent unf}
Let $f:\GS{n}{p}$ be an $\eqA$-finite germ. We say that an origin-preserving one-parameter unfolding $F(x,t)=\big(f_t(x),t\big)$ is \textit{good} if there exists a representative $F:U\to W\times T$, where $U$ is an open neighbourhood of $S\times\{0\}$ in $\CC^n\times\CC$ and $W,$ $T$ are open neighbourhoods of the origin in $\CC^{p},$ $\CC$ respectively, such that
		\begin{enumerate}[\itshape(i)]
			\item $F$ is finite,
			\item $f_t^{-1}(0)=S$, for all $t\in T$,
			\item \label{i3:excellentunf}$f_t$ is locally stable on $W-\left\{0\right\}$, for all $t\in T$.
		\end{enumerate}
		If, in addition, $f_t$ has no 0-stable singularities on $W-\{0\}$ (i.e., stable singularities whose isosingular locus is 0-dimensional, see \cite[Definition 5.2]{Houston2010} or \cite[p. 93]{Mond2020}) 
      we say that $f_t$ is \textit{excellent} \textit{(in Gaffney's sense)}.
\end{definition}

\begin{figure}[ht]
	\centering
		\includegraphics[width=1.00\textwidth]{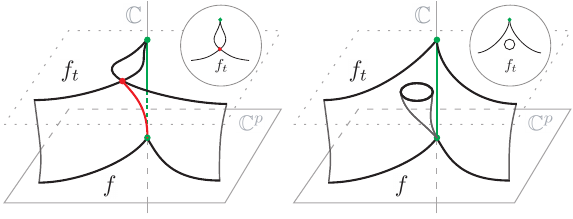}
		\caption{A non-excellent unfolding (left) due to the presence of a 1-dimensional stratum, of stable or unstable points, distinct from the parameter axis (red and green, respectively), and an excellent unfolding (right) with only one stratum of dimension one (green). }
%
	\label{fig:Excellent unfolding}
\end{figure}

Excellent unfoldings play an important role in the theory of equisingularity of families of germs (see, for example,  \cite{Gaffney1993}). See \cref{fig:Excellent unfolding} for a representation of an excellent unfolding.


Now we can state Houston's conjecture on excellent unfoldings (see \cite[Conjecture 6.2]{Houston2010}). 
\begin{conjecture}[Houston's conjecture]
Let $f:\GS{n}{p}$, $n<p$, be $\eqA$-finite of corank $1$ and  let $F(x,t)=\big(f_t(x),t\big)$ be an origin-preserving one-parameter unfolding. Then, if $\widetilde{f_t}$ are the germs $f_t:\GS{n}{p}$, $\mu_I(\widetilde{f_t})$ constant implies $F$ excellent.
\end{conjecture}

We give a counterexample to this conjecture in \cref{ex:disgustingsplitting}. We also give the correct modification of the conjecture, which we also prove, in \cref{hc}. 

\begin{example}\label{ex:disgustingsplitting}
The map germ $f:(\CC^7,0)\to(\CC^{11},0)$ with one-parameter unfolding given by $$ f_t(\underline{x},y)=\big(\underline{x},\ x_6y+y^6,\ x_5y+y^5,\ x_4y+y^4,\ x_3^2y+y^3,\ x_2y+y^2(x_1^2-tx_1)\big)$$
has $D^2(f_t)$ with constant Milnor number (at the origin) for all $t$ and 
$$D^3(f_t)=\left\{\left(\underline{0}_6;0,0,0\right)\right\} \cup \left\{\left(t,\underline{0}_5;0,0,0\right)\right\}\subset \CC^{6+3},$$
so $f_t$ has two instabilities when $t\neq0$ and $\mu_I(f_t;0)>0$ is constant (see \cref{thm: conustedesmui}).
\end{example}

\begin{theorem}[cf. {\cite[Theorem 4.3]{GimenezConejero2022}}]\label{hc}
Let $f:\GS{n}{p}$, $n<p$, be $\eqA$-finite of corank $1$ and  let $F(x,t)=\big(f_t(x),t\big)$ be an origin-preserving one-parameter unfolding. 
Assume that $f_t$ has neither 
 strongly contractible instabilities nor unexpected homology. Then, if $\widetilde{f_t}$ are the germs $f_t:\GS{n}{p}$, $\mu_I(\widetilde{f_t})$ constant implies $F$ excellent.
\end{theorem}
\begin{proof}
We can assume that $f$ is unstable, otherwise it is trivial.
By \cite[Corollary 5.9]{Houston2010}, we only need to prove goodness of the unfolding and that, for the induced germ from $f_t$ at zero (not necessarily equal to $\widetilde{f_t}$), either $s(f_t;0)\leq d(f_t;0)$ for all $t$ or both are constant.

That the unfolding $F$ is good follows from the hypotheses. Indeed, we have conservation of the $\mu_I$ because there is no unexpected homology (see \cref{thm: conservation muinui}) and every instability has positive $\mu_I$, by the absence of strongly contractible instabilities. Hence, the constancy of the image Milnor number implies that there is only one instability (the origin) for each $f_t$, this confirms \cref{i3:excellentunf} of \cref{def: excellent unf}. 


If $s(f_t;0)>d(f_t;0)$ for some $t$, then $d(f_t)=\kappa=\left\lfloor \frac{p}{p-n}\right\rfloor$ by \cref{rem: s>d} and $$\mu_{d(f)+1}^\Alt(f_t;0)>0.$$ Because we are in conditions of conservation of $\mu_I$, $\mu_I$ and every $\mu_k^\Alt$ is 
 upper semi-continuous; so necessarily $$\mu_{d(f)+1}^\Alt(f;0)>0$$ and $s(f)>d(f)=\kappa$. We only need to prove that $s(f_t;0)=s(f)$. As $F$ is origin-preserving, we already know that $s(f_t;0)\geq s\big(\widetilde{f_t}\big)=s(f)$ but, if the inequality is strict, then
$$ \mu_{d(f)+1}^\Alt(f_t;0)>\mu_{d(f)+1}^\Alt(f), $$
which is absurd by upper semi-continuity of $\mu_I$.
\end{proof}

It is not difficult to see the reason to add the hypothesis on strongly contractible instabilities: once we know that there are instabilities with $\mu_I=0$, we could find perturbations $f_t$ with strongly contractible instabilities splitting from the instability at the origin (see a non-trivial example in \cref{ex:disgustingsplitting}). 

The condition on unexpected homology is also easy to understand: there could be constancy of the $\mu_I$ while having instabilities away from the origin because of the creation of unexpected homology (that makes instabilities with $s(\bullet)>d(\bullet)$, possibly at the origin, by \cref{thm: conservation muinui}). In other words, without strongly contractible instabilities and unexpected homology, the presence of one instability away from the origin implies that $\mu_I$ cannot be constant at the origin; but, with unexpected homology, it could happen that the reduction of the $\mu_I$ at the origin is counterbalanced with more branches because of the creation of unexpected homology.

\begin{corollary}
Let $f:\GS{n}{p}$, $n<p$, be $\eqA$-finite of corank $1$ and  let $F(x,t)=\big(f_t(x),t\big)$ be an origin-preserving one-parameter unfolding. Assume that $D^k(f)=\varnothing$ for $d_k<0$. Then, if $\widetilde{f_t}$ are the germs $f_t:\GS{n}{p}$, $\mu_I(\widetilde{f_t})$ constant implies $F$ excellent.
\end{corollary}
\begin{proof}
In this case, there are no strongly contractible instabilities nor unexpected homology in any perturbation (see \cref{rem:disg,rem: unexpected}).
\end{proof}

We can ask about the converse of this result, i.e., if excellency implies the constancy of the image Milnor number. This was addressed in \cite[Proposition 4.4]{GimenezConejero2022}, where the author and Nuño-Ballesteros provide a proof for the case $(n,p)=(1,2)$ and $(n,p)=(2,3)$. There is, however, a small mistake in the proof of the case $(n,p)=(2,3)$: one needs to assume that $f_t$ has an instability, as the following example shows.

\begin{figure}[htb]
	\centering
		\includegraphics[width=0.60\textwidth]{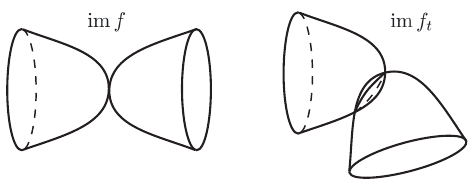}
	\caption{Deformation $f_t$ that is excellent but does not have constant image Milnor number.}
	\label{fig:Tangency2}
\end{figure}

\begin{example}\label{ex:doubletwo}
The tangent double point germ given in \cref{ex: singularidad libre}, $f:\big(\CC^2,\{0,0'\}\big)\to(\CC^3,0)$ given by 
$$  \left\{\begin{aligned} 
(x,y)&\mapsto (x,y,x^2+y^2)\\
(x',y')&\mapsto(x',y',0)
\end{aligned}\right. $$
has an excellent one-parameter unfolding such that $f_t$ is stable at every point and which also preserves the origin (see \cref{fig:Tangency2}):
$$ \left\{\begin{aligned} 
(x,y)&\mapsto \big(x,y\cos(t)-\sin(t)(x^2+y^2),\cos(t)(x^2+y^2)-\sin(t)y\big)\\
(x',y')&\mapsto(x',y',0)
\end{aligned}\right. $$
This implies that there are excellent families $f_t$ with $\mu_I(f_t;0)$ not constant. However, the \textsl{tangent double point} is the only example of this phenomenon in these dimensions, any other excellent family has an instability.
\end{example}

The same authors presented \cite[Conjecture 4.6]{GimenezConejero2022} for the remaining cases, which we prove now false by a counterexample:

\begin{example}\label{ex:counterextangent}
A \textit{fold} singularity is an $\eqA$-finite germ $f:(\CC^n,0)\to(\CC^p,0)$ 
 that is $\eqA$-equivalent to $f(\underline{x},y)=\big(\underline{x},y^2,yp_{n+1}(\underline{x},y),\dots,yp_{p}(\underline{x},y)\big)$, where $p_i(\underline{x},y)$ are polynomials. 

If we consider a fold germ $f:\Gzero{3}{4}$, none of its perturbations has zero-stable points (stable singularities whose isosingular locus is 0-dimensional). This comes from the fact that zero-stable points in corank one are in one-to-one correspondence with multiple point spaces $D^k(\bullet)^\sigma$ such that $d_k^\sigma=0$. However, these germs only have double points, so $D^{2+\ell}(f)=\varnothing$ for $\ell>0$, and $d_2,d_2^{(1\; 2)}\geq1$.

One can certainly find many excellent families, for example given by stable perturbations $f_s$ that fix the origin (which can be assumed), such that $\mu_I(f)>\mu_I(f_s)=0$.
\end{example}

\cref{ex:counterextangent} works for any pair $(n,n+1)\neq(1,2),(2,3)$, so \cite[Proposition 4.4]{GimenezConejero2022} is, indeed, optimal if we add the detail about the presence of an instability in every $f_t$.

\begin{note}
This kind of counterexamples with fold maps was found, independently but also simultaneously, by the author, Oset Sinha and Nuño-Ballesteros after fruitful conversations between them. 
\end{note}

Regarding the general question in dimensions $(n,p)$, we can give a complete answer.

%

\begin{proposition}
In the context of  $\eqA$-finite map germs $f:\GS{n}{p}$ of corank one, for $n<p<2n$ there are excellent families with non-constant image Milnor number. Furthermore, for $n=2p$ excellency implies constancy of $\mu_I$.
\end{proposition}
\begin{proof}
Observe that $p=2n$ is equivalent to $d_2=0$, and $p<2n$ to $d_2>0$. Hence, we can divide by cases depending on $d_2$.


If $d_2>1$, one can always find fold germs that are unstable and they all give counterexamples 
 (as we have seen in \cref{ex:counterextangent}). 

If $d_2=0,1$, the presence of strongly contractible instabilities is not a problem: they do not change the image Milnor number. The case $d_2=0$ is easy, $\mu_I$ only depends on the number of points in the fiber of the zero-dimensional \textsc{icis} $D^2(f)$ (see \cref{conjrgc}), so a change in $\mu_I$ implies presence of more than one instability or some zero-dimensional stable points, making the unfolding not excellent.

If $d_2=1$, then $p=2n-1$. We take the idea of \cref{ex:doubletwo}. Hence, we consider the germ $f:\big(\CC^n,\{0,0'\}\big)\to(\CC^{2n-1},0)$ given by the generalization of the tangent double point given in \cref{ex: singularidad libre,ex:doubletwo}, i.e.,
$$  \left\{\begin{aligned} 
(x_1,\dots,x_n)&\mapsto (x_1,\dots,x_{n},x_1,\dots,x_{n-2},x_{n-1}^2+x_{n}^2)\\
(x_1',\dots,x_n')&\mapsto(x_1',\dots,x_{n}',0,\dots,0)
\end{aligned}\right. $$
The equations of $D^2(f)\subset\CC^{n+n}$ are $x_i=x_i'$, $x_1=\dots=x_{n-2}=0$ and $x_{n-1}^2+x_{n}^2=0$, which gives an \textsc{icis} of dimension one. Furthermore, $D^k(f)=\varnothing$ for $k>2$, so $f$ is $\eqA$-finite by the Marar-Mond criterion \cref{thm: Marar-Mond crit}. Since there are no points fixed by the permutation $(1\: 2)\in\Sigma_2$ (as in the tangent double point \cref{ex: singularidad libre}), i.e., $D^2(f)^{(1\:2)}=\varnothing$, any stable deformation on one parameter $f_s$ is excellent because there are no 0-stable singularities (similarly to what was explained in \cref{ex:doubletwo}). 
 
\end{proof}

\section{Conclusions and open problems}\label{sec: conjectures}
As we have seen, the technique leading to \cref{conjrgc,thm: conustedesmui} allows us to translate the invariants $\mu_I$ and $\nu_I$ into Milnor numbers of \textsc{icis} in corank one. In corank two or higher one could do something similar, but Euler-Poincaré characteristics of some varieties would need to be used. Indeed, as we have pointed out in \cref{rem: marar-mond formulae}, this is the correct generalization of the Marar-Mond formulae for $\mu_I$, see \cite{Marar1989}.

Several applications 
 of this technique were shown, such as the proof Houston's conjecture \cref{hc}. The main goal of this paper is being able to bypass the difficulty of working with $\mu_I$ and use the easier-to-handle Milnor numbers of \textsc{icis}. Furthermore, we have clarified many traits of deformations of $\eqA$-finite germs $f:(\CC^n,S)\to(\CC^p,0)$, with $n<p$, using these ideas. However, we have only seen the corank one case. It seems natural to conjecture the following, generalizing the notion of strongly contractible germs to any corank.

\begin{oproblem}
Show that there are unstable $\eqA$-finite germs $f:(\CC^n,0)\to(\CC^p,0)$, $n<p$, of corank bigger than one such that have contractible stable perturbation.
\end{oproblem}

Moreover, it is unclear what are the dimensions where we can find strongly contractible \textsl{multi-germs} (even of corank one). We know, by \cref{ex:multiSS}, that \cref{thm:disgusting?} is not sharp for multi-germs.  There are some results regarding this: by \cite[Theorem 2.15]{GimenezConejero2023a} (corank one case was also given in \cite[Theorem 3.9]{GimenezConejero2022}), if $p=n+1$ then every stable perturbation of an unstable germ $f$ has non-trivial homology, i.e., $\mu_I(f)>0$. 

It is, however, reasonable to think that there are no strongly contractible germs if, and only if, $d_\kappa=0$ (i.e., if $\frac{p}{p-n}\in\ZZ^+$) because it seems that strongly contractible germs and unexpected homology (controlled by by \cref{cor: sing iff alt}) have something in common: see \cref{ex:multiSS} where a strongly contractible bi-germ given from $f_t$ at the origin appears together with some unexpected homology of a deformation of another germ $h$ (recall the discussion on \cref{subs:unexhom}).  

\begin{oproblem}
There are unstable $\eqA$-finite germs $f:(\CC^n,S)\to(\CC^p,0)$, $n<p$, such that $\im f_s$ is contractible 
 if, and only if, $\frac{p}{p-n}\notin\ZZ^+$; for $f_s$ the stable perturbation.
\end{oproblem}

Also, the weakened version of this conjecture for corank one is the generalization of \cref{thm:disgusting?}: 

\begin{oproblem}
There are strongly contractible multi-germs $f:(\CC^n,S)\to(\CC^p,0)$ of corank one, $n<p$, if, and only if, $\frac{p}{p-n}\notin\ZZ^+$.
\end{oproblem}

Going back to the corank one case, and regarding unexpected homology, it is not clear that mono-germs have deformations with unexpected homology. Indeed, every example, such as \cref{ex:multiSS}, seems to need more than one branch to find a deformation with unexpected homology.

\begin{oproblem}
There are no $\eqA$-finite mono-germs $f:\Gzero{n}{p}$, with $n<p$, such that they have a deformation with unexpected homology. 
\end{oproblem}

%

This would eliminate the hypothesis of unexpected homology of \cref{hc} for the case of mono-germs and, indeed, \cref{ex:multiSS} is not a counterexample since it is a bi-germ.
\newline


Finally, the initial motivation of this paper was solving the following conjecture.

\begin{conjecture}\label{conj: splitting}
For $\eqA$-finite germs $f:(\CC^n,S)\to(\CC^{n+1},0)$, there cannot be coalescence of instabilities in one-parameter deformations. In other words, a one-parameter deformation $f_t$ of $f$ cannot have more than one instability and contractible image, $\im f_t$, at the same time.
\end{conjecture}

Regarding coalescence of instabilities in the general case $n<p$, it is not clear that there are mono-germs with coalescence of strongly contractible germs when \cref{thm:disgusting?} holds with equality: $d_\kappa-\kappa=0$. The germ given in \cref{ex:disgustingsplitting}, which has this coalescence, has $d_\kappa-\kappa=1$.

\begin{oproblem}\label{oprob: coal strong contr}
There are no $\eqA$-finite germs $f:(\CC^n,0)\to(\CC^p,0)$ of corank one such that they have coalescence of strongly contractible instabilities and $d_\kappa-\kappa=0$.
\end{oproblem}

\begin{remark}
The Mond conjecture is a long-standing open problem of map germs $f:(\CC^n,S)\to(\CC^{n+1},0)$, stating that the $\eqA_e$-codimension is lower or equal that the $\mu_I$, with equality in the quasi-homogeneous case. The literature on this problem is vast and points to a positive answer, but the general statement remains open (see \cite{deJong1991,Mond1995,Houston1998,Houston1999a,Cooper2002,Sharland2014,Sharland2019} among many others). It seems reasonable to think that there is a version of \cref{thm: conustedesmui} for the $\eqA_e$-codimension, translating it into invariants of \textsc{icis}. This would lead to a possible proof of the conjecture in corank one, based on the known inequality $\mu\geq\tau$ for \textsc{icis} of positive dimension (see \cite{Greuel1975,Greuel1980,Looijenga1985} and the nice survey \cite{Greuel2019}). In turn, this could lead to a full proof with some ideas from \cite{Bobadilla2019}. However, it is not clear what is the analogue 
 of \cref{thm: conustedesmui} for the $\eqA_e$-codimension.

Furthermore, \cref{conj: splitting} has deep connections with this conjecture as both relate the same invariants (see \cite[Section 7.3]{Robertothesis}): one-parameter unfoldings correspond to curves in the parameter space of a versal unfolding. 
\end{remark}

\bibliographystyle{abbrv}
\bibliography{C:/Users/rgc19/Documents/Textos-Matematicas/Bibliografias/bibtex/bib/mybibs/FullBib.bib} 
\end{document}